\documentclass[12pt]{elsarticle}
\usepackage{graphicx}
\usepackage{epsfig}
\usepackage{amsmath}
\usepackage{amssymb}
\usepackage{amsfonts}
\usepackage{euscript}
\begin{document}
\bibliographystyle{plain}
\newcommand{\bea}{\begin{eqnarray}}
\newcommand{\eea}{\end{eqnarray}}
\newcommand{\bfm}[1]{{\mbox{\boldmath{$#1$}}}}
\newcommand{\bfmN}{{\mbox{\boldmath{$N$}}}}
\newcommand{\bfmx}{{\mbox{\boldmath{$x$}}}}
\newcommand{\bfmv}{{\mbox{\boldmath{$v$}}}}
\newcommand{\se}{\setcounter{equation}{0}}
\newtheorem{corollary}{Corollary}[section]
\newtheorem{proposition}{Proposition}[section]
\newtheorem{example}{Example}[section]
\newtheorem{definition}{Definition}[section]
\newtheorem{theorem}{Theorem}[section]
\newtheorem{lemma}{Lemma}[section]
\newtheorem{remark}{Remark}[section]
\newtheorem{result}{Result}[section]
\newcommand{\vtwo}{\vskip 4ex}
\newcommand{\vthree}{\vskip 6ex}
\newcommand{\vfour}{\vspace*{8ex}}
\newcommand{\hone}{\mbox{\hspace{1em}}}
\newcommand{\hon}{\mbox{\hspace{1em}}}
\newcommand{\htwo}{\mbox{\hspace{2em}}}
\newcommand{\hthree}{\mbox{\hspace{3em}}}
\newcommand{\hfour}{\mbox{\hspace{4em}}}
\newcommand{\von}{\vskip 1ex}
\newcommand{\vone}{\vskip 2ex}
\newcommand{\n}{\mathfrak{n} }
\newcommand{\m}{\mathfrak{m} }
\newcommand{\q}{\mathfrak{q} }
\newcommand{\aF}{\mathfrak{a} }

\newcommand{\kl}{\mathcal{K}}
\newcommand{\p}{\mathcal{P}}
\newcommand{\Lt}{\mathcal{L}}
\newcommand{\bv}{{\mbox{\boldmath{$v$}}}}
\newcommand{\bc}{{\mbox{\boldmath{$c$}}}}
\newcommand{\bx}{{\mbox{\boldmath{$x$}}}}
\newcommand{\br}{{\mbox{\boldmath{$r$}}}}
\newcommand{\bs}{{\mbox{\boldmath{$s$}}}}
\newcommand{\bb}{{\mbox{\boldmath{$b$}}}}
\newcommand{\ba}{{\mbox{\boldmath{$a$}}}}
\newcommand{\bn}{{\mbox{\boldmath{$n$}}}}
\newcommand{\bp}{{\mbox{\boldmath{$p$}}}}
\newcommand{\by}{{\mbox{\boldmath{$y$}}}}
\newcommand{\bz}{{\mbox{\boldmath{$z$}}}}
\newcommand{\be}{{\mbox{\boldmath{$e$}}}}

\newcommand{\bP}{{\mbox{\boldmath{$P$}}}}

\newcommand{\M}{\mathcal{M}}
\newcommand{\R}{\mathbb{R}}
\newcommand{\Q}{\mathbb{Q}}
\newcommand{\Z}{\mathbb{Z}}
\newcommand{\N}{\mathbb{N}}
\newcommand{\C}{\mathbb{C}}
\newcommand{\xar}{\longrightarrow}
\newcommand{\ov}{\overline}
 \newcommand{\rt}{\rightarrow}
 \newcommand{\om}{\omega}
 \newcommand{\wh}{\widehat }
 \newcommand{\wt}{\widetilde }
 \newcommand{\g}{\Gamma}
 \newcommand{\lm}{\lambda}

\newcommand{\eN}{\EuScript{N}}
\newcommand{\ncom}{\newcommand}
\newcommand{\norm}{\|\;\;\|}
\newcommand{\inp}[2]{\langle{#1},\,{#2} \rangle}
\newcommand{\nrm}[1]{\parallel {#1} \parallel}
\newcommand{\nrms}[1]{\parallel {#1} \parallel^2}
\title{Adaptive Pad\'e-Chebyshev Type Approximation to Piecewise Smooth Functions}
\author{ S. Akansha\footnote{akansha@math.iitb.ac.in} and S. Baskar\footnote{baskar@math.iitb.ac.in}\\
Department of Mathematics, \\Indian Institute of Technology Bombay,\\ Powai, Mumbai - 400076. India.}
\maketitle{}
\begin{center}{\bf Abstract}\end{center}
The aim of this article is to study the role of piecewise implementation of Pad\'e-Chebyshev type approximation in minimising Gibbs phenomena in approximating piecewise smooth functions. A piecewise Pad\'e-Chebyshev type (PiPCT) algorithm is proposed and an $L^1$-error estimate for at most continuous functions is obtained using a decay property of the Chebyshev coefficients.  An advantage of the PiPCT approximation is that we do not need to have an {\it a prior} knowledge of the positions and the types of singularities present in the function.  Further, an adaptive piecewise Pad\'e-Chebyshev type (APiPCT) algorithm is proposed in order to get the essential accuracy with a relatively lower computational cost. Numerical experiments are performed to validate the algorithms. The numerical results are also found to be well in agreement with the theoretical results. Comparison results of the PiPCT approximation with the singular Pad\'e-Chebyshev and the robust Pad\'e-Chebyshev methods are also presented.

\section{Introduction}
%

Often in applications, we come across the problem of approximating a non-smooth function, for instance, while approximating the entorpy solution of a hyperbolic conservation law that involves shocks and rarefactions, and the compression of an image that involves edges. The main challenge in such an approximation problem is the occurrence of the well known Gibbs phenomena in the approximant near a singularity of the target function.  The methods that are proven to have higher order of convergence for smooth target functions tend to have lower order of convergence near jump discontinuities (Gottlieb and Shu \cite{got-shu_97a}).

Nonlinear approximation methods are often preferred to linear approximation procedures in order to achieve convergence in a vicinity of singularities with a higher order of convergence (see DeVore \cite{dev_98a}).  Multiresolution  schemes have been developed based on the Harten's ENO procedure, where an appropriate set of nodes are chosen locally to minimize the effect of the singularity on the approximation (Arandiga {\it et al.} \cite{ara-etal_05a}).  Other approach is to construct the approximation using a standard linear approximation and then use filters and mollifiers to reduce the Gibbs phenomena (see Tadmor \cite{tad_07a}).  By knowing some more information about the singularities of the target function, say the location and/or the nature of the singularity, one can develop efficient and accurate methods (see for instance, Gottlieb and Shu  \cite{got-shu_95a}, Driscoll and Fornberg \cite{dri-for_01a} and Lipman and Levin \cite{lip-lev_10a}).

Another nonlinear approach for approximating a piecewise smooth function is the rational approximation, in particular, the Pad\'e approximation.  The main advantage of rationalising a truncated series is to achieve a faster convergence in the case of approximating analytic functions  (see Baker {\it et al.} \cite{bak_65a, bak-etal_96a}) when compare to a polynomial approximation.  The advantage of rational approximation also lies in its property of auto error correction (see Litvinov \cite{lit_03a}), in which the error in all the intermediate steps is compensated at the final step of rationalisation. This is particularly useful for approximating piecewise smooth functions.  Geer \cite{gee_95a} used a truncated Fourier series expansion of a periodic even or odd piecewise smooth function to construct a Pad\'e approximant, which has been further developed for a general function by Min {\it et al.} \cite{min-etal_07a}.  Hesthaven {\it et al.} \cite{hes-etal_06a} proposed a Pad\'e-Legendre interpolation method to approximate a piecewise smooth function.  Kaber and Maday \cite{kab-mad_05a} studied  the convergence rate of a sequence of Pad\'e-Chebyshev approximation for sign-function.  These nonlinear methods minimise Gibbs phenomena significantly, however the convergence at the jump discontinuity is rather slow.
In order to accelerate the convergence rate further, Driscoll and Fornberg \cite{dri-for_01a} (also see Tampos {\it et al.} \cite{tam-etal_12a}) developed a singular Pad\'e approximation to a piecewise smooth function using finite series expansion coefficients or function values. By knowing the jump locations, this method captures the discontinuities of a piecewise smooth function very accurately. The methods mentioned above are global approximations.  But in many applications (like in constructing numerical schemes for partial differential equations), it is desirable to have a local approximation technique that minimises the Gibbs phenomena without the explicit knowledge of the type, magnitude, and location of the singularities of the target function.

%
%

In this article, we are mainly concerned to have a local approximation technique based on the Pad\'e-Chebychev type approximation that is suitable for approximating a piecewise smooth function. To this end, we propose a piecewise implementation of the Pad\'e-Chebychev type approximation (denoted by PiPCT approximation) and demonstrate numerically that the proposed method captures the isolated singularities (including jump discontinuities) of a non-smooth function.  For a given partition $P_N$  of an interval $[a,b]$ consisting of $N$ subintervals, and the degrees of the numerator and the denominator polynomials in each subinterval as $N$-dimensional vectors $\bfm{n}_p$ and $\bfm{n}_q$, respectively, the  PiPCT algorithm computes the Pad\'e-Chebyshev approximation (of order $[n_p^j/n_q^j]$) of the truncated Chebyshev series of a function $f\in L^1([a,b])$ in each subinterval $I_j$, $j=0,1,\ldots, N-1$, of the given partition with Chebyshev coefficients being approximated with $n$ quadrature points.  
We prove $L^1$-convergence of the PiPCT approximation as $N\rightarrow \infty$ under the assumption
 that in each subinterval, the function $f^{(k-1)}$, for a given positive integer $k$, is absolutely continuous and $f^{(k)}$ is bounded.
Further, we demonstrate numerically the convergence of the PiPCT approximation in a vicinity of jump discontinuity with a higher order of convergence. The key advantage of the PiPCT approximation is that it reduces Gibbs phenomena significantly without specifying the location of the jump discontinuity a priori. We compare the results of PiPCT with the singular Pad\'e approximant \cite{dri-for_01a} and the Robust Pad\'e-Chebyshev approximant \cite{gon-etal_11a, gon-etal_13a} (both computed on the whole interval, referred as global approximation) and found that the PiPCT approximation captures the singularities sharply, which is in comparison with the singular Pad\'e approximation.

Our numerical experiments show that the proposed PiPCT algorithm captures the jump discontinuities very sharply without the knowledge of the location of the jump.  Also, our theoretical study shows that we do not need a piecewise implementation of PCT approximation in the regions where the target function is more smoother.  This motivates us to look for an adaptive algorithm which gives an optimal partition that is required to approximate the target function with accuracy as much as we obtain in the piecewise algorithm with uniform partition. This obviously needs the location of the singularities of the target function.  Many work had been devoted in the literature to develop singularity indicator, see for instance, Banerjee and Geer \cite{ban-gee_98a}, Barkhudaryan {\it et al.} \cite{bar-etal_07a}, Eckhoff \cite{eck_95a}, Kvernadze \cite{kve_04a}.  Using the result that the genuine poles appear near the singularities (Baker {\it et al.} \cite{bak-etal_61a}), we identify the intervals where the denominator polynomial of the PCT approximation almost vanishes and define these intervals as {\it bad-cells}. We further refine the partition in an isotropic greedy manner (for details see chapter 3 of \cite{dev-kun_09a}).   
 This results in an adaptive piecewise Pad\'e-Chebyshev type (APiPCT) algorithm, where the singularity regions are identified with a less additional computational cost. 

The article is organised as follows: decay property of coefficients of the Chebyshev series expansion of a function is discussed in section \ref{Cheb.Approx.sec}. The PiPCT algorithm is developed in section \ref{PCT.approximation.sec}, where a  $L^1$-error estimate and order of accuracy results for PiPCT approximant are also discussed. Numerical experiments to validate the PiPCT algorithm are present in section \ref{numerical.comparison.sec}, where we compare the results of PiPCT algorithm with singular Pad\'e and robust Pad\'e-Chebyshev methods. Section \ref{APiPCT.sec} is devoted towards the proposed PiPCT based adaptive algorithm (APiPCT). Numerical evidence on the performance of APiPCT is presented in section \ref{numerical.experiment.APiPCT.sec}.   We also implemented the robust Pad\'e-Chebyshev approximation in our adaptive algorithm (denoted by APiRPCT) and compared the results with APiPCT in this section.  Finally, in section \ref{comment.FD.sec}, we discuss numerically the presence of the Froissart doublets in the PiPCT approximation and their role in the accuracy of the approximant.

\section{Chebyschev Approximation}\label{Cheb.Approx.sec}

For a given function $f \in L_\omega^2[-1,1],$ with $\omega(x)=1/\sqrt{1-x^2},$ the Chebyshev series representation of $f$ is given by
\begin{equation}\label{CMap.eq}
f(x)=\sideset{}{'}\sum_{k=0}^{\infty} c_k T_k(x), ~~ x\in [-1,1],
\end{equation}
where the prime in the summation indicates that the first term is halved, $T_k(x)$ denotes the Chebyshev polynomial of degree $k,$ and $c_k$, $k=0,1,\ldots$, are the Chebyshev coefficients given by 
\bea\label{CSRepCoeff.eq}
c_k = \frac{2}{\pi}\langle f, T_k\rangle_w.
\eea
Here, $\langle\cdot,\cdot\rangle_ \omega$ denotes the weighted $L^2$-inner product.

The Chebyshev coefficients \eqref{CSRepCoeff.eq} are approximated using the Gauss Chebyshev quadrature formula (see \cite{mas-han_03a,riv_74a}) given by
\begin{equation}\label{eq:approxcoeff}
c_{k,n}:=\frac{2}{n}\sum_{l=1}^{n}f(t_l)T_k(t_l) ,~~k=0,1,\ldots,
\end{equation}
where the quadrature points $t_l$, $l=1,2,\ldots, n$ are the Chebyshev points given by
\bea\label{ChebPts.eq}
t_l	= \cos\left(\dfrac{\left(l+0.5\right) \pi}{n}\right), \hspace{1cm} l=1, 2, \ldots,n. 
\eea
We use the notation
$$ \mathsf{C}_{d,n}[f](x):=\sideset{}{'}\sum_{k=0}^{d} c_{k,n} T_k(x)$$
for the truncated Chebyshev series up to degree $d$ with approximated coefficients involving $n$ quadrature points and the series is denoted by
$$ \mathsf{C}_{\infty,n}[f](x):=\sideset{}{'}\sum_{k=0}^{\infty} c_{k,n} T_k(x).$$

The Chebyshev series representation of any function $f \in L_\omega^2[a,b]$ can be obtained using a change of variable
\begin{equation}\label{eq:bijection}
	\mathsf{G}(y) = a+(b -a)\frac{(y+1)}{2}, ~~ y\in [-1,1].
\end{equation}
A $L^1$-error estimate  in approximating $f$ by $\mathsf{C}_{d,n}[f]$ can be obtained from the following decay estimates of the Chebyshev coefficients:
\begin{theorem}[Majidian \cite{maj_17a}]\label{thm:coeffbound}
For some integer $k\ge 0$, let $f^{(k-1)}$ be an absolutely continuous function on the interval $[a,b]$ and let $f^{(k)}$ be of bounded variation. 
If $\|f^{(k)}\|_T = V_k< \infty,$ where
\begin{equation}\label{eq:chebnormscale}
\|f\|_T := \int_{0}^{\pi} \left|{f'\left(\mathsf{G}(\cos \theta )\right)}\right|d\theta,
\end{equation}
then the following inequalities hold:
\begin{enumerate}
\item If $k = 2s$ for some integer $s\geq 0$ ($k$ is even)
\begin{equation}\label{eq:coeffbound1}
|c_n| \leq \left(\dfrac{b-a}{2}\right)^{2s+1}  \dfrac{2V_k}{\pi\displaystyle \prod_{j = -s}^{s}(n + 2j)},~~n\ge k+1.
\end{equation}
\item If $k = 2s+1$ for some integer $s\geq 0$ ($k$ is odd)
\begin{equation}\label{eq:coeffbound2}
|c_n| \leq \left(\dfrac{b-a}{2}\right)^{2s+2}  \dfrac{2V_k}{\pi\displaystyle \prod_{j = -s}^{s+1}(n + 2j-1)},~~n\ge k+1.
\end{equation}
\end{enumerate}
\end{theorem}
Note that, the estimates in the above theorem are not providing any interesting information for $f$ being analytic and in that case we have the following well known decay estimate (see Rivlin \cite{riv_74a} and also see Xiang {\it et al.} \cite{xia-etal_10a}):
\begin{theorem}\label{thm:AnalyticCoffbound}
Let $f: [a,b] \rightarrow \mathbb{R}$ be a real analytic function and let $f$ has an  analytic extension in the ellipse $C_\rho$ with foci $\pm 1$ and the sum of major and minor axes equals $\rho >1$. If $|f(z)| \leq C$ for $z \in C_\rho,$ then for each $j \geq 0$
\begin{equation*}
|c_j| < \dfrac{2C}{\rho^j}.
\end{equation*}
\end{theorem}
The required error estimate can now be proved using the above decay estimates of the Chebyshev coefficients:
\begin{theorem}\label{thm:errestimate}
Under the hypothesis of Theorem \ref{thm:coeffbound} and for $n-1\ge k\ge 1,$ we have
$$\|f - \mathsf{C}_{d,n}[f]\|_1
 \leq C_{d,n},$$
 where
\begin{enumerate}
\item for $d=n+l$, $-n\le l \le 0$,
\begin{equation}\label{eq:totalL1err}
C_{d,n}= \left\{
\begin{array}{@{}l@{\thinspace}l}
\left(\dfrac{b-a}{2}\right)^{2s+2} \dfrac{4V_k}{k \pi } \left(\dfrac{1}{\displaystyle \prod_{i = -s}^{s-1}(n+l+2i+1)} + \dfrac{1}{\displaystyle \prod_{i = -s+1}^{s}(n+l+2i+1)} \right),  & \text{ if } k = 2s\medskip\\		
\left(\dfrac{b-a}{2}\right)^{2s+3} \dfrac{4V_k}{k\pi} \left(\dfrac{1}{\displaystyle \prod_{i = -s}^{s}(n+l+2i)} + \dfrac{1}{\displaystyle \prod_{i = -s}^{s}(n+l+2i+1)} \right), & \text{ if } k = 2s+1, \\
\end{array}
\right.
\end{equation}
\item and for $d = n+l, l > 0$
\begin{equation}
C_{d,n} = \left\{
\begin{array}{@{}l@{\thinspace}l}
\left(\dfrac{b-a}{2}\right)^{2s+2} \dfrac{6V_k}{\pi k } \left(\dfrac{1}{\displaystyle \prod_{i = -s}^{s-1}(n-l+2i)} + \dfrac{1}{\displaystyle \prod_{i = -s+1}^{s}(n-l+2i)} \right),  & \text{ if } k = 2s\medskip\\		
\left(\dfrac{b-a}{2}\right)^{2s+3} \dfrac{6V_k}{k\pi} \left(\dfrac{1}{\displaystyle \prod_{i = -s}^{s}(n-l+2i-1)} + \dfrac{1}{\displaystyle \prod_{i = -s}^{s}(n-l+2i)} \right). & \text{ if } k = 2s+1. \\
\end{array}
\right.
\end{equation}
\end{enumerate}
\end{theorem}
\textbf{Proof:} 
Using the well-known result (see Rivlin \cite{riv_74a}) 
$$c_d - c_{d,n} = \sum_{j=1}^\infty (-1)^j\big(c_{2jn-d} + c_{2jn+d}\big),$$
we have
\begin{align}\label{Est01.eq}
\|f - \mathsf{C}_{d,n}[f]\|_1
& \leq (b-a)\left[\sum_{j = d+1}^{\infty} \left|c_j\right| + \sum_{j = 1}^{\infty} \sum_{k = 2jn -d}^{2jn+d}\left|c_k\right|\right].
\end{align}
For $d = n+l, -n\leq l \le 0$, we have
\begin{align*}
\sum_{j = n+l+1}^{\infty} \left|c_j\right| + \sum_{j = 1}^{\infty} \sum_{k = 2jn - (n+l)}^{2jn+(n+l)}\left|c_k\right|
& \leq 2\sum_{j = n+l+1}^{\infty} \left|c_j\right|.
\end{align*}
For  $d = n+l, l > 0$, the right hand side of \eqref{Est01.eq} can be rewritten as
\begin{align*}
\sum_{j = n+l+1}^{\infty} \left|c_j\right| + \sum_{j = 1}^{\infty} \sum_{k = 2jn - (n+l)}^{2jn+(n+l)}\left|c_k\right| 
& = 2\sum_{j = n+l+1}^{\infty} \left|c_j\right| + \sum_{j = 1}^{\infty} \sum_{k = (2j-1)n-l}^{(2j-1)n+l}\left|c_k\right|.
\end{align*}
Using the decay estimate from Therorem \ref{thm:coeffbound} and then using the telescopic property of the resulting series (see also Majidian \cite{maj_17a}), we can arrive at the required estimates.
\qed
\begin{remark}
From the above theorem, we see that for a fixed $n$ (as in the hypothesis), the upper bound $C_{d,n}$ decreases for $d=n-l$ and increases for $d=n+l$ as $l\in [0,n]$ increases, and attains its minimum at $d=n$, which gives the interpolating polynomial at Chebyshev nodes.  Further, we see that $C_{n-l-1,n} = \frac{3}{2}C_{n+l,n}$,
however computationally  $\mathsf{C}_{n-l-1,n}[f]$ is more efficient than $\mathsf{C}_{n+l,n}[f].$ A similar situation occurs in the case of  Pad\'e-Chebyshev type approximation as shown in subsection \ref{higherPCT.lowerPCT.ssec} and numerically shown in  \ref{degree.adaptation.subs}.\qed
\end{remark}
On similar lines, the error estimate in the case when $f$ is analytic can be proved using Theorem \ref{thm:AnalyticCoffbound} (also see Xiang {\it et al.} \cite{xia-etal_10a}).  
\begin{theorem}
Under the hypotheses of Theorem \ref{thm:AnalyticCoffbound},
\begin{enumerate}
\item  for each $n \geq 0$ and $d = n-l, 0 \leq l \leq n$, we have
\begin{equation}
\|f - \mathsf{C}_{d,n}[f]\|_1 \leq \dfrac{4C(b-a)}{(\rho-1)\rho^{n-l}},
\end{equation}
\item and for each $n \geq 0$ and $d = n+l, l > 0$, we have
\begin{equation}
\|f - \mathsf{C}_{d,n}[f]\|_1 \leq \dfrac{6C(b-a)}{(\rho-1)\rho^{n-l-1}}.
\end{equation}
\end{enumerate}
\end{theorem}
Also note that the hypotheses of the above theorems restrict us to use the error estimates only for functions $f$ that does not involve a jump discontinuity.  For a discontinuous function, the Chebyshev approximant may develop Gibb's phenomena in a vicinity of the jump discontinuity.  The Gibb's phenomena can be considerably reduced (although not fully removed) if we go for a rationalization of the Chebyshev approximant, for instance, the Pad\'e-Chebyshev approximation.  
\section{Pad\'e-Chebyshev type approximation}\label{PCT.approximation.sec}
In this section, we recall the basic construction of the Pad\'e-Chebyshev type approximation of a given function $f\in L^2_w[-1,1]$ and further propose a piecewise implementation of this approximation.

For $x\in [-1,1],$ we use the notation 
\bea\label{CComplexNotation.eq}
\mathsf{C}_{\infty}[f](z)
:= \sideset{}{'}\sum_{k=0}^{\infty} c_k z^k;~~z = e^{i\cos^{-1}(x)},
\eea
which is a complex power series defined on the unit circle centered at the origin in the complex plane whose coefficients are real and are given by \eqref{CSRepCoeff.eq}.  It can be seen that the real part of this series is precisely the Chebyshev series given by \eqref{CMap.eq}.  

 For given non-negative integers $n_p\ge n_q \geq 1$, a rational function 
 \bea\label{PadeApprox.eq}
 \mathsf{R}_{n_p,n_q}(z) := \dfrac{P_{n_p}(z)}{Q_{n_q}(z)}
 \eea
  with numerator polynomial $P_{n_p}(z)$ of degree  $\leq n_p$  and denominator polynomial $Q_{n_q}(z)$ of degree $\leq n_q$ with $Q_{n_q} \neq 0$ satisfying (see \cite{bec-mat_15a,dri-for_01a,tam-etal_12a})
\begin{equation} \label{eq:P_linear}
Q_{n_q}(z)\mathsf{C}_{\infty}[f](z)-P_{n_p}(z)= \mathit{O}(z^{n_p+n_q+1}), \hspace*{.5cm} z\rightarrow 0,
\end{equation}
is called a Pad\'e approximant of $\mathsf{C}_{\infty}[f](z)$ of order $[n_p/n_q]$. Such a Pad\'e approximation exists and the real part of the $ \mathsf{R}_{n_p,n_q}(z)$ is an approximation of $f(x)$, which is referred as a {\it Pad\'e-Chebyshev type} (PCT) approximant of $f$ (see \cite{tam-etal_12a}).

The coefficients of the polynomial $Q_{n_q}(z)$ are obtained as a solution of the Toeplitz system \cite{dri-for_01a,hes-etal_06a,tam-etal_12a,bec-mat_15a}
\begin{equation*}
\begin{bmatrix}
c_{n_p+1} & c_{n_p} &\cdots & c_{n_p-n_q+1} \\
c_{n_p+2} & c_{n_p+1} &\cdots & c_{n_p-n_q+2} \\
\vdots &\vdots&\vdots &\vdots \\
c_{n_p+n_q} &c_{n_p+n_q-1} & \cdots & c_{n_p} 
\end{bmatrix} \begin{bmatrix}
q_0 \\ q_1 \\ \vdots \\ q_{n_q}
\end{bmatrix}
= \begin{bmatrix}
0 \\ 0 \\ \vdots \\ 0
\end{bmatrix}.
\end{equation*}
In the matrix notation, we write the above system as
\begin{equation}\label{eq:linearsystem}
A_{n_p,n_q}\bfm{q} = \bfm{0},
\end{equation}
where
$$A_{n_p,n_q}=(c_{i-k}),~~i=n_p+1,n_p+2\ldots,n_p+n_q,~~k=0,1,\ldots,n_q.$$
Once the coefficients of the denominator polynomial are known, the coefficients of the numerator polynomial  $P_{n_p}(z)$ can be computed using the matrix vector multiplication
\begin{equation}
\begin{bmatrix}
p_0 \\ p_1 \\ \vdots\\ p_{n_q}\\ \vdots \\ p_{n_p}
\end{bmatrix} =
\begin{bmatrix}
c_{0}/2 & 0 & \cdots & 0 \\
c_{1} &c_0/2 & \cdots & 0 \\
\vdots & & & \\
c_{n_q} & c_{n_q-1} & \cdots & c_{0}/2 \\
\vdots & & & \\
c_{n_p} & c_{n_p-1} & \cdots & c_{n_p-n_q} 
\end{bmatrix} \begin{bmatrix}
q_0 \\ q_1 \\ \vdots \\ q_{n_q}
\end{bmatrix} \label{eq:PCnum}.
\end{equation}
A unique PCT can be obtained if $A_{n_p,n_q}$ is of full rank.
A Pad\'e-Chebyshev type approximation of $f$ of order $[n_p/n_q]$ can be computed for the given set of Chebyshev coefficients $\{c_0,c_1,\ldots,c_{n_p+n_q}\}$.  Since these coefficients often cannot be obtained exactly, we use the approximated coefficients $\{c_{0,n},c_{1,n},\ldots,c_{n_p+n_q,n}\}$ obtained by the Gauss-Chebyshev quadrature formula \eqref{eq:approxcoeff} to compute PCT and denote it by $ \mathsf{R}_{n_p,n_q}^{n}.$
\subsection{Computing higher order PCT approximants using lower order PCT}\label{higherPCT.lowerPCT.ssec}
This subsection presents an interesting property of $\mathsf{R}_{n_p,n_q}^{n}$ (also see Cuyt and Wuytack \cite{cuy-wuy_87a} for a similar result) which motivated us in choosing a suitable degree of the numerator polynomial in subsection \ref{degree.adaptation.subs}.
\begin{proposition}\label{prop:1}
	Let $f\in L^{2}_{\omega}[a,b]$ and consider the approximated Chebyshev series $f(t) \approx \sum_{k=0}^{\infty}c_{k,n}T_k(t)$ for $t\in [-1,1].$ Let $\mathsf{R}_{n_p,n_q}^{n}(z)$ be the unique Pad\'e approximation of $\mathsf{C}_{\infty,n}[f](z)$ of order $[n_p/n_q]$. Then
for Pad\'e approximations of order $[(n-j-1)/n_q]$ and $[(n+j)/n_q]$, $j = 0,1,\ldots,n-n_q-1$, the corresponding denominator coefficient vectors $\bfm{q}_{n_q}^{n-j-1}$ and $\bfm{q}_{n_q}^{n+j}$ satisfies
		\begin{equation}\label{eq:Denrelationgen}
		\bfm{q}_{n_q}^{n-j-1} = R_{n_q+1}\bfm{q}_{n_q}^{n+j},
		\end{equation}
where $R_{n_q+1}$ is the 'flip' matrix of size $(n_q+1) \times (n_q+1)$ with $1$ on anti-diagonal and $0$ everywhere else.
\end{proposition}
We use the following anti-symmetric property of the approximated Chebyshev series coefficients to prove Theorem \ref{prop:1}.
\begin{lemma}\label{lem:2}
	The approximated Chebyshev series coefficients \eqref{eq:approxcoeff} of $f$ using $n$ quadrature points satisfy, for $k=1,3,5,\ldots$,
	\begin{equation}\label{eq:result1}
	\left.\begin{array}{rlr}
	c_{kn+j,n}  &= -c_{kn-j,n}, & \hspace{.3cm} j = 1,2,\ldots,2n, \\
	c_{kn,n} &= 0.&
	\end{array}\right\}
	\end{equation}
\end{lemma}
\textbf{Proof: } It is clear to see that $c_{kn,n} = 0,$ for $k=1,3,5,\ldots.$ Now multiplying both sides by $2T_j(t_l)$ and using the trigonometric identity 
\begin{equation*}
2\cos A \cos B = \cos (A+B) +\cos(A-B),
\end{equation*}
we obtain the required result.\qed

\noindent\textbf{Proof of Theorem \ref{prop:1}: }
By the construction of the Pad\'e-Chebyshev approximants of order $[(n-j-1)/n_q]$ and $[n+j/n_q]$, for $j = 0,1,\ldots,n-n_q-1$, the denominator coefficient vectors $\bfm{q}_{n_q}^{n-j-1}$ and $\bfm{q}_{n_q}^{n+j}$ can be computed by solving the linear systems
\begin{equation}\label{eq:linearsys(n-j-1)}
A_{n-j-1,n_q}\bfm{q}_{n_q}^{n-j-1} = \bfm{0} \hspace{1cm} \mbox{and} \hspace{1cm} 
A_{n+j,n_q}\bfm{q}_{n_q}^{n+j} = \bfm{0},
\end{equation} 
respectively.  Using Lemma \ref{lem:2}, we can write
$$A_{n+j,n_q} = (-c_{i+k}), \hspace{1cm} \mbox{for } i = n-j-1,n-j-2,\ldots,n-j-n_q,~k=0,1,\ldots,n_q.$$
A direct observation shows that  by flipping rows and columns of the matrix $A_{n-j-1,n_q}$ and then multiplying by -1 we obtain the right hand side matrix in the above representation.  Thus, we obtain
\begin{equation}\label{eq:prop1}
	A_{n+j,n_q} = -R_{n_q}A_{n-j-1,n_q}R_{n_q+1},
\end{equation}
where $R_{n_q}$ and $R_{n_q+1}$ are `flip' matrices of size $n_q$ and $n_q+1$, respectively.  Let $x \in \ker (A_{n+j,n_q})$ then from \eqref{eq:prop1}, $R_{n_q}A_{n-j-1,n_q}R_{n_q+1}x = 0$ or $A_{n-j-1,n_q}R_{n_q+1}x  = 0.$	Therefore for any vector $x \in \ker(A_{n+j,n_q})$ implies $R_{n_q+1}x \in \ker(A_{n-j-1,n_q})$ or vise versa.	This completes the proof.\qed
\begin{corollary}
In particular, if $j=0$ in Theorem \ref{prop:1}, the \eqref{eq:Denrelationgen} can be written as
\begin{equation}\label{eq:Denrelation}
		\bfm{q}_{n_q}^{n-1} =\pm \bfm{q}_{n_q}^{n}.
\end{equation}
\end{corollary}
\noindent\textbf{Proof:}
The proof follows from the fact that the vector $\bfm{q}_{n_q}^{n}$ becomes an eigenvector of the flip matrix $R_{n_q+1}$.\qed
\subsection{Piecewise Implementation}\label{piecewise.implementation.ssec}

Consider a function $f\in L^1[a,b]$. Let us first discretize the interval $I:=[a,b]$ into $N$ cells,  denoted by $I_j:=[a_j,b_j]$, $j=0,1,\ldots,N-1$, where
$$a=a_0<b_0=a_1<b_1=a_2<\cdots<b_{N-2}=a_{N-1}<b_{N-1}=a_N=b$$
and denote the partition as $P_N:=\{a_0,a_1,\ldots,a_N\}.$

For the given integers $n$ and $N$, and $(N-1)$-tiples $\bfm{n_p}=(n_{p}^{0},\ldots,n_{p}^{N-1})$ and $\bfm{n_{q}}=(n_{q}^{0},\ldots,n_{q}^{N-1}),$   
construct the PiPCT approximation of $f$ as follows:
\begin{enumerate}
\item Generate $n$ Chebyshev points $\{t_l~:~ l=1,2,\cdots,n\},$ given by \eqref{ChebPts.eq}, in the reference interval $[-1,1]$. \\
\item For each $j=0,1, \ldots, N-1$, consider the bijection map $\mathsf{G}_j : [-1,1] \rightarrow I_j$ given by
	\begin{equation}\label{eq:bijection}
	\mathsf{G}_j(y) = a_j+(b_j -a_j)\frac{(y+1)}{2}.
	\end{equation}
Obtain the approximate Chebyshev coefficients in the $j^{\rm th}$ cell, denoted by $c_{k,n}^j,$ for $k=0,1,\cdots,n_{p}^{j}+n_{q}^{j},$ using the Gauss Chebyshev quadrature formula \eqref{eq:approxcoeff} with the values of $f$ evaluated at $\mathsf{G}_j(t_l)$, $l=1,2,\cdots,n$.\\
%
%
\item Define the {\it piecewise Pad\'e-Chebyshev type}  (PiPCT) approximation of $f$ in the interval $I$ with respect to the given partition as
\bea\label{PPC.eq}
\mathsf{R}^{n,N}_{n_p,n_q}(x):=
\begin{cases}
\mathsf{R}_{n_{p}^0,n_{q}^0}^{n}(x),&\quad\text{if }x \in [a_0,b_0),\\
\mathsf{R}_{n_{p}^1,n_{q}^1}^{n}(x) ,&\quad\text{if } x\in [a_1,b_1),\\
\vdots \\
\mathsf{R}_{n_{p}^{N-1},n_{q}^{N-1}}^{n}(x),&\quad\text{if } x \in [a_{N-1},b_{N-1}],
\end{cases}
\eea
where  $\mathsf{R}_{n_{p}^j,n_{q}^j}^{n}(x),$ for $x\in I_j$, denotes the PCT approximant of $f|_{I_j}$, $j=0,1,\ldots,N-1$.
\end{enumerate}
Using Theorem \ref{thm:errestimate} we can obtain an error bound for the PiPCT approximation. 
\begin{theorem}\label{errorestimate.thm}
 Assume the hypothesis of Theorem \ref{thm:coeffbound} for $f|_{I_j}$, $j=0,1,\ldots, N-1$. Then, for each fixed integer $n-1 \geq k\ge 1$ and for $n_p^j\geq n_q^j\geq 1$,  $j=0,1,\ldots, N-1$, the $L_1$-order of convergence of PiPCT approximation is at least $k+1$ as $N \rightarrow \infty.$
\end{theorem}
\textbf{Proof: } 
We have
\bea
\big\|f-\mathsf{R}^{n,N}_{n_p,n_q}\big\|_1 &\le& \sum_{j=0}^{N-1}\big\|f|_{I_j}-\mathsf{R}_{n_{p}^j,n_{q}^j}^{n}\big\|_1 \nonumber\\
&\le& \sum_{j=0}^{N-1}\big\|f|_{I_j}-\mathsf{C}_{d^j,n}[f|_{I_j}]\big\|_1 +  \sum_{j=0}^{N-1}\big\| \mathsf{C}_{d^j,n}[f|_{I_j}]-\mathsf{R}_{n_{p}^j,n_{q}^j}^{n}\big\|_1,\nonumber
\eea
where $d^j = n_{p}^j+n_{q}^j.$ From the construction of the Pad\'e-Chebyshev type approximant we see that the second term in the above inequality vanishes and hence we have
\bea
\big\|f-\mathsf{R}^{n,N}_{n_p,n_q}\big\|_1 
&\le& \sum_{j=0}^{N-1}\big\|f|_{I_j}-\mathsf{C}_{d^j,n}[f|_{I_j}]\big\|_1.
\eea
Without loss of generality, let us take $d=d^j$, for $j=0,1,\ldots, N-1$. 
Now using $L_1$-error bound from Theorem \ref{thm:errestimate} for truncated Chebyshev series approximation we have
\begin{enumerate}
	\item for $d = n+l, -n\leq l\leq 0$
{\small	\begin{equation}\label{eq:totalL1errPade1}
\|f - \mathsf{R}^{n,N}_{n_p,n_q}\|_1 \leq \left\{
	\begin{array}{@{}l@{\thinspace}l}
	\left(\dfrac{h}{2}\right)^{2s+1} \dfrac{4V_k(b-a)}{k \pi } \left(\dfrac{1}{\displaystyle \prod_{i = -s}^{s-1}(n+l+2i+1)} + \dfrac{1}{\displaystyle \prod_{i = -s+1}^{s}(n+l+2i+1)} \right),  & \text{ if } k = 2s\medskip\\		
\left(\dfrac{h}{2}\right)^{2s+2} \dfrac{4V_k(b-a)}{k\pi} \left(\dfrac{1}{\displaystyle \prod_{i = -s}^{s}(n+l+2i)} + \dfrac{1}{\displaystyle \prod_{i = -s}^{s}(n+l+2i+1)} \right), & \text{ if } k = 2s+1. \\
	\end{array}
	\right.
	\end{equation}}
	\item and for $d = n+l, l > 0,$
	\begin{equation}\label{eq:totalL1errPade2}
\|f - \mathsf{R}^{n,N}_{n_p,n_q}\|_1 \leq \left\{
	\begin{array}{@{}l@{\thinspace}l}
	\left(\dfrac{h}{2}\right)^{2s+1} \dfrac{6V_k(b-a)}{\pi k } \left(\dfrac{1}{\displaystyle \prod_{i = -s}^{s-1}(n-l+2i)} + \dfrac{1}{\displaystyle \prod_{i = -s+1}^{s}(n-l+2i)} \right),  & \text{ if } k = 2s\\		
\left(\dfrac{h}{2}\right)^{2s+2} \dfrac{6V_k(b-a)}{k\pi} \left(\dfrac{1}{\displaystyle \prod_{i = -s}^{s}(n-l+2i-1)} + \dfrac{1}{\displaystyle \prod_{i = -s}^{s}(n-l+2i)} \right), & \text{ if } k = 2s+1, \\
	\end{array}
	\right.
	\end{equation}
\end{enumerate}
where $h = \dfrac{b-a}{N}$.
\qed

Note that PiPCT involves two parameters, namely, the smoothness parameter $k$ and the discretization parameter $h$. From the above estimates, we observe that the upper bound of the $L^1-$error in PiPCT  is smaller for a smoother function, and it tends to zero as $h\rightarrow 0$. Therefore we can conclude that PiPCT algorithm is well suited for approximating a sufficiently smooth function.  However, we do not have a theoretical justification for the accuracy of the algorithm in a vicinities of singularities of a piecewise smooth function.  
\section{Numerical Comparison}\label{numerical.comparison.sec}
There are two tasks to be addressed numerically.  One is to examine the performance of the PiPCT approximation proposed in the previous section and the other one is to study numeically  the convergence rate of the proposed algorithm in order to show that the rate of convergence obtained theoretically in Theorem \ref{errorestimate.thm} is achieved numerically.

As a first task, we give numerical evidence that the proposed algorithm captures singularities (of all order) of a function without a visible Gibbs phenomena.   
\begin{example}{\rm 
Consider the piecewise smooth function
	\begin{equation}\label{discondeg.eq}
	f(x) = 
	\begin{cases}
	x^3, &\quad\text{if } x \in [-1,-0.4),\\
	x^2+1, &\quad\text{if } x \in [-0.4,0.4),\\
	1.16 - (x-0.4)^{\frac{1}{2}}, &\quad\text{if } x \in [0.4,1].
	\end{cases}
	\end{equation}
The function $f$ involves a jump discontinuity at $x=-0.4,$ whereas at $x=0.4$ the function is continuous but not differentiable (referred in this article as {\it point singularity}) as shown in Figure \ref{fig:PiPChe_discondeg}(a). We fix the number of quadrature points as $n=200$ so that the error due to the quadrature formula is considerably reduced.  

The purpose of this numerical experiment is to study the accuracy of the PiPCT approximant in terms of its parameters.   The purpose is also to compare the performance of PiPCT with recently proposed global Pad\'e-based algorithms for approximating functions with singularities, namely, the singular
Pad\'e-Chebyshev (SPC) (see Driscoll and Fornberg \cite{dri-for_01a}, and Tampos {\it et al.} \cite{tam-etal_12a}) and the robust Pad\'e-Chebyshev type algorithms (RPCT) (see Gonnet \cite{gon-etal_11a, gon-etal_13a}) and altogether with global PCT approximant.

The proposed PiPCT algorithm is performed to approximate the function $f$ by taking the number of subintervals $N = 512$ and then fixing $n = 200$ and $n_p^j=n_q^j=20$, for $j=0,1,\ldots,N-1,$ in all the uniformly discretised subintervals of a partition. For all global algorithms we fixed the parameters as $n_p = 20 = n_q$, $n = 512 \times 200$, and $N = 1$.  
We can see from Figure \ref{fig:PiPChe_discondeg}(a) that in the smooth region these approximants are well in agreement with the exact function but in the vicinity of the singularities $x=-0.4$ and $x=0.4$, SPC is able to capture the jump discontinuity but not point singularity (see the zoomed boxes).  Also, it is clear form the figure that the RPCT and the global PCT approximants perform almost the same in a small vicinity of the singularities.  From this figure, we can observe a significant role of piecewise implementation of the PCT approximation in capturing the singularities without knowing the location of the singularities {\it a priori} unlike in the case of SPC.

\noindent {\it Numerical discussion as $N$ varies:} Figure \ref{fig:PiPChe_discondeg}(b) depicts the peaks of the pointwise error (as explained in Driscoll and Fornberg \cite{dri-for_01a}) for $N=2^k,$ $k=0,1,3,5,7,8,$ and $9.$ Figure \ref{fig:PiPChe_discondeg}(c) depicts the maximum error in the vicinity of $x=-0.4$ (in $-\!\!\!-\!\!\!$o$\!\!\!-\!\!\!-$ symbol) and in the vicinity of $x=0.4$ (in $-\!\!\!-\!\!\!*\!\!\!-\!\!\!-$ symbol) with logarithmic scale in the $y$-axis.
\begin{figure}[t]
	\centering
\includegraphics[height=14.cm,width=14cm]{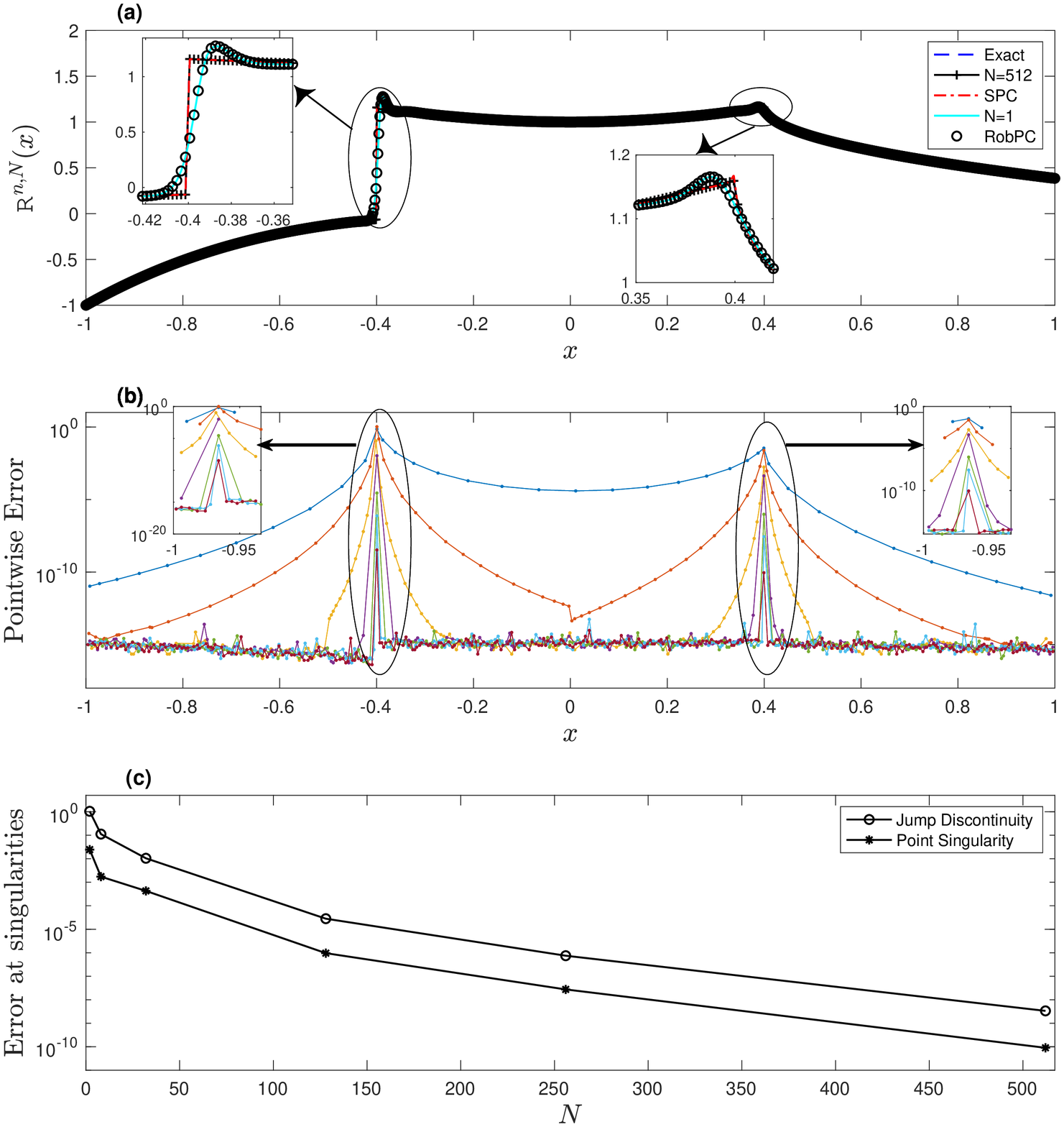}\vspace{-0.3in}
\caption{(a) Depicts comparison of PiPCT approximant of $f$ given by  \eqref{discondeg.eq} with the global Pad\'e-Chebyshev based algorithms, (b) depicts the peaks of pointwise error in approximating $f$ by PiPCT for $N = 2^k, k = 0,1,3,5,7,8,9$, and (c) depicts the $L^\infty$-error  in approximating $f$ by PiPCT  in a vicinity of $x = -0.4$ and $x = 0.4$.}
\label{fig:PiPChe_discondeg}
\end{figure}
\vspace{0.1in}

Note that $\mathsf{R}^{102400,1}(x)$ is a global approximation whereas the later one $\mathsf{R}^{200,512}(x)$ is a piecewise approximation.   In both the cases we have given the values of the function at 102400 points.  In the case of the global approximant, these points are the Chebyshev points in the interval $[-1,1]$ and we can clearly observe an oscillation near the discontinuity as shown in the zoomed box of Figure \ref{fig:PiPChe_discondeg}(a) in the vicinity of $x=-0.4$.  The same kind of behaviour is also observed in the vicinity of the point singularity at $x=0.4$.  Note that the error estimates obtained in Theorem \ref{errorestimate.thm} indeed shows that the sequence of piecewise approximants converges for functions that are at least continuous.  Though we do not have a theorem that gives convergence for functions involving jump discontinuities, the piecewise approximant in this example captures the singularities accurately including the jump discontinuity at $x=-0.4$. Moreover, it is evident from Figure \ref{fig:PiPChe_discondeg}(b) and (c) that the sequence of piecewise approximants (in this example) tends to converge to the exact function as $N\rightarrow \infty$. 

Table \ref{ErrorNVaries.tab} compares the $L^1$-error and the numerical order of accuracy of the piecewise Chebyshev and the PiPCT approximants in the interval $[0.2,1]$ where the function $f$ has a point singularity at $x=0.4$.  Here, we can see the obvious advantage of using Pad\'e-Chebyshev type approximation when compared to the Chebyshev approximation. Also, we observe that the numerical order of convergence is well in agreement with the theoretical result given in Theorem \ref{errorestimate.thm} with $k=1,$ $d=40$, and $n=200$.
\begin{table}
\caption{The $L^1$-error and the numerical order of accuracy of the piecewise Chebyshev and the PiPCT approximants of the function \eqref{discondeg.eq} in the interval $[0.2,1]$ as $N$ varies. We have taken $n_p=n_q=20$ in all the subintervals and $n=200.$}\vspace{0.1in}
\centering
\begin{tabular}{|l| c|c| c| c| c|  r|}
\hline
\multicolumn{1}{|c|}{$N$}&\multicolumn{2}{|c|}{Piecewise Chebyshev}&\multicolumn{2}{|c|}{Piecewise Pad\'e-Chebyshev type}\\
\hline
&$L^1$-error & order & $L^1$-error & order\\ 
\hline
2  &  0.603860 & --     &0.032616 &--     \\
\hline
8  &  0.03022378662560210039 & 2.725905     &0.00064588620006190815 &3.569908     \\
\hline
32  &  0.00109697967811300638 & 2.552230     &0.00002635315776778789 &2.462154     \\
\hline
128  &  0.00003324170968713331 & 2.564731     &0.00000001505864286582 &5.477419     \\
\hline
256  &  0.00000323501556305392 & 3.380088     &0.00000000021392558412 &6.171918     \\
\hline
512  &  0.00000004032820231387 & 6.343662     &0.00000000000035272088 &9.270412     \\
\hline
\end{tabular}
\label{ErrorNVaries.tab}
\end{table}

\begin{figure}[t]
	\centering
\includegraphics[height=13.cm,width=12cm]{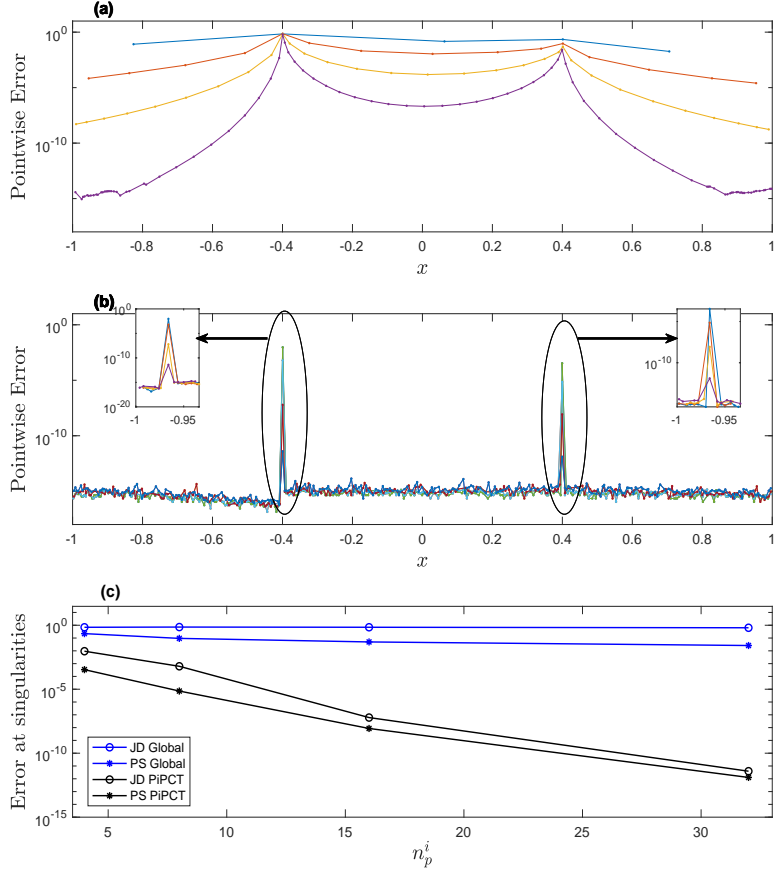}
\caption{Comparison between local and global PC algorithms for $N = 512,$ $n = 200$ and $n_p = n_q = 2,4,6,8,16,32$. (a) depicts the peaks of the pointwise error $|f(x) - \mathsf{R}^{N\times n,1}_{n_p,n_q}(x)|$, (b) depicts the peaks of the pointwise error $|f(x) - \mathsf{R}^{n,N}_{n_p,n_q}(x)|$, and (c) depicts the convergence of both the algorithms in a vicinity of $ x = -0.4$ and $x = 0.4$.}
\label{fig:PiPCheandChe_discondeg}
\end{figure}

\noindent {\it Numerical discussion as $n_p (=n_q)$ varies:}
Let us take $n_p^j=n_q^j=n_p$, $j=0,1,\ldots,N-1,$ and study the numerical convergence of the sequence of PiPCT approximantions $\mathsf{R}^{n,N}_{n_p,n_p}$ of the function $f,$ given by \eqref{discondeg.eq}, as $n_p$ varies. In order to have a clear advantage of using the piecewise approximation, 
we also study numerically the convergence of the global Pad\'e-Chebyshev type approximantions $\mathsf{R}^{n\times N,1}_{n_p,n_p}$.

The peaks of the pointwise errors $|f(x)-\mathsf{R}^{n\times N,1}_{n_p,n_p}(x)|$ and $|f(x)-\mathsf{R}^{n,N}_{n_p,n_p}(x)|$ are depicted in Figures \ref{fig:PiPCheandChe_discondeg} (a) and (b), respectively. Here we have taken $n=200$, $N=512$ (therefore, in both the cases, we use the function values at 102400 grid points) and varied $n_p = 2, 4, 8, 16$ and $32.$ We clearly observe that the pointwise error in the global implementation of the Pad\'e-Chebyshev type approximation does not seems to be converging at the jump discontinuity (at $x=-0.4$) as $n_p$ increases and the convergence is slow at the point singularity at $x=0.4$. This behaviour is more apparent in Figure \ref{fig:PiPCheandChe_discondeg} (c). On the other hand, from Figure \ref{fig:PiPCheandChe_discondeg} (b) and (c), we see that the convergence in the piecewise approximation is significantly faster both at the jump discontinuity and at the point singularity as $n_p$ increases.
%
\qed
}
\end{example}
 In the above example we have seen that the numerical order of convergence is well in agreement with the theoretical result given in \eqref{eq:totalL1errPade1} when $k=1$ ({\it i.e.} when $k$ is odd).  
 
 In the following example, we show that this result also holds numerically for $k=2$.
 \begin{example}\label{SmoothPPC.ex}{\rm
 Consider the $C^1$ function
 $f(x) = x|x|$ 
 for $x\in [-1,1]$.  We use both the piecewise Chebyshev and the PiPCT approximations. The $L^1$-error and the order of convergence are tabulated in Table \ref{SmoothPPC.tab}. The results are well in agreement with the theoretical results of Theorem \ref{errorestimate.thm} for $k=2$.
 }
 \qed
\begin{table}
\caption{The $L^1$-error and the numerical order of accuracy of the piecewise Chebyshev and the PiPCT approximants for the function given in Example \ref{SmoothPPC.ex}, as $N$ varies. We have taken $n_p=n_q=2$ in all the subintervals and $n=200.$}\vspace{0.1in}
\centering
\begin{tabular}{|l| c|c| c| c| c|  r|}
\hline
\multicolumn{1}{|c|}{$N$}&\multicolumn{2}{|c|}{Piecewise Chebyshev}&\multicolumn{2}{|c|}{Piecewise Pad\'e-Chebyshev}\\
\hline
&$L^1$-error & order & $L^1$-error & order\\ 
\hline
2  &  0.00000000000001916544 & --     &0.00000000000002741904 &--     \\
\hline
4  &  0.00000000000000282010 & 3.751452     &0.00000000000000335724 &4.111221     \\
\hline
8  &  0.00000000000000028910 & 3.875159     &0.00000000000000031289 &4.037213     \\
\hline
16  &  0.00000000000000003397 & 3.366981     &0.00000000000000003508 &3.440576     \\
\hline
\end{tabular}
\label{SmoothPPC.tab}
\end{table}
 \end{example}
\section{An Adaptive Algorithm}\label{APiPCT.sec}
In section \ref{numerical.comparison.sec}, we demonstrated numerically the performance of the PiPCT approximation in capturing singularities of a function accurately.
 However, the numerical results depicted in Figure \ref{fig:PiPChe_discondeg} suggests that we need a sufficiently finer discretization to get a good accuracy in a vicinity of the singularities.  Such a finer discretization is not needed in the regions where the function is relatively smoother.  This can also be observed in the error estimates  \eqref{eq:totalL1errPade1}-\eqref{eq:totalL1errPade2}. This motivates us to look for an adaptive implementation of the PiPCT algorithm suggested in subsection \ref{piecewise.implementation.ssec}.

The main idea of our adaptive algorithm is to identify the subintervals (referred as {\it badcells}) where the function has a singularity and then bisect the subinterval.  Therefore, as a first step, we need a strategy to identify the badcells.
\subsection{Singularity Indicator}

In general rational (Pad\'e) approximations are better than (see for instance, \cite{kab-mad_05a,bak-pet_82a}) the polynomial approximations in the case of approximating non-smooth functions. A natural concern about  Pad\'e approximations is the poles of the approximant. It may happen that the approximant has poles at places where the function has no singularities. Such poles are often referred as the spurious poles  (for precise definition see \cite{sta_98a}).  In addition to such spurious poles, the Pad\'e-Chebyshev type approximants develop poles in a sufficiently small neighborhood of a singularity of the function (Baker {\it et al.} \cite{bak-etal_61a}). Poles are intractable on a computer because of the presence of round-off error. Nevertheless, we have observed numerically the following result:
\begin{result}\label{res:jumplocate}
Let $f\in L^1[-1,1]$, which is at most continuous, with an isolated singularity at $x_0 \in [-1,1]$. For a given  $\epsilon>0$, there exists an integer $n_{p_0}>0$,  such that
\begin{equation}\label{Adaptive.Condition.eq}
|Q_{n_{p}}(z)|<\epsilon,~~\mbox{\rm for all }~~n_p\ge n_{p_0},
\end{equation}
whenever $|x_0-{\rm Re}(z)|<\delta,$ for some $\delta=\delta(n_p,\epsilon)>0.$ 
Here, $Q_{n_p}(z)$ $($for all $z\in \mathbb{C}$ such that $|z|=1)$ is the denominator polynomial in the Pad\'e-Chebyshev type approximation $\mathsf{R}_{n_p,n_p}$ of $f$.
\end{result}
We illustrate this result in the following numerical example.
\begin{figure}[t]
	\centering
\includegraphics[height=5.cm,width=15cm]{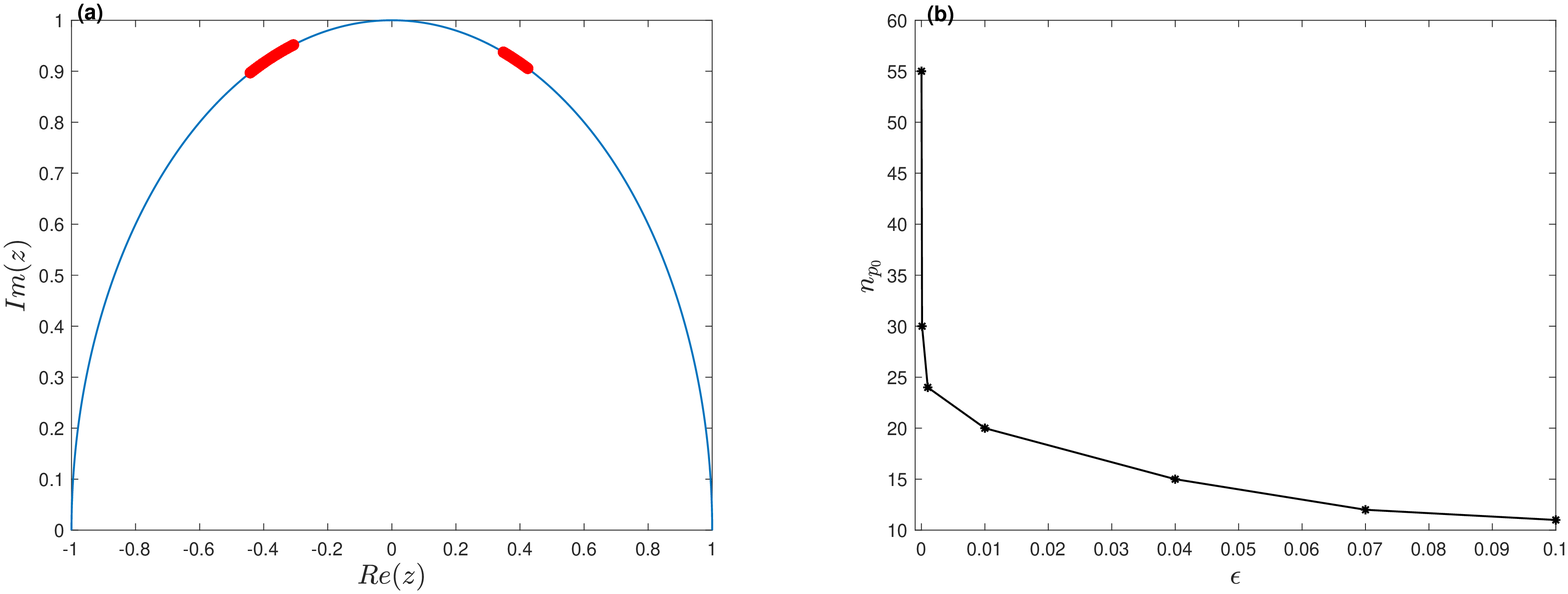}
\caption{(a) depicts the points (in red `o' symbol) on the unit circle which satisfies $|Q_{20}(z)| < 10^{-2}$, and (b) depicts the $\epsilon$ values and the corresponding $n_{p_0}$ which satisfies \eqref{Adaptive.Condition.eq}.}
\label{fig:Singularity.Location.discondeg}
\end{figure}
\begin{example}
{\rm 
Consider the piecewise smooth function $f$ given by \eqref{discondeg.eq}, where the function has a jump discontinuity at $x=-0.4$ and a point singularity at $x=0.4.$ Figure \ref{fig:Singularity.Location.discondeg} (a) depicts the points (in red `o' symbol) at which $|Q_{20}(z)|\le 10^{-2}.$ We observe that the real part of these points are accumulated in neighborhoods of the points $x=-0.4$ and $x=0.4$.   Further, Figure \ref{fig:Singularity.Location.discondeg} (b) depicts the graph whose $x$-coordinate represents different values of $\epsilon$ and the $y$-coordinate is taken to be the corresponding values of $n_{p_0}$, which are the minimum values of the degree of the denominator polynomial $Q_{n_p}$, for which the condition \eqref{Adaptive.Condition.eq} holds numerically.
\qed
}
\end{example}
Based on the above numerical observation, we define a notion of  {\it badcells} in a given partition as follows:
\begin{definition}\label{eps.bad.cell.def}{\rm
Let $f\in L^1[a,b]$ and let $P_N$ be a given partition. For a given $\epsilon>0$ and an integer $m>0$, a subinterval $I_j=[a_j,b_j]$ with $a_j,b_j\in P_N$ is said to be an {\it $\epsilon$-badcell} if  $|Q_{m}(z)|<\epsilon$ for some $z\in \mathbb{C}$ on the unit circle with ${\rm Re}(z)\in \mathsf{G}_j^{-1}(I_j)$.
}
\end{definition}
\subsection{Generation of an Adaptive Partition}\label{generation.adaptive.partition.ssec}
In this subsection, we propose an algorithm to generate a partition which consists of finer discretisation in a vicinity of singularities.

Consider a function $f\in L^1[a,b]$.  Choose an $\epsilon>0$ sufficiently small, a tolerance parameter $\tau>0$, a positive integers $n$ sufficiently large, and $m(\ll n/2)$ sufficiently small with
 $\bfm{n}_p=\bfm{n}_q=(m,m).$ 

Let us take $N_0=2$, $a_0=a$ and $b_0=b.$ Let $I^{(0)}_0$ and $I^{(0)}_1$ be the subintervals of the partition $P_{N_0}:=\{a_0,(a_0+b_0)/2,b_0\}.$ 

Perform the following steps for $j=0,1,2,\cdots$:
\begin{enumerate}
\item Obtain the denominator polynomials, denoted by $Q_{m,k}(z),$ $k=0,1$, of the Pad\'e-Chebyshev type approximant of the function $f$ for the interval $I^{(j)}_k$.
\item Check if the subinterval $I^{(j)}_k$ is an $\epsilon$-badcell using Definition \ref{eps.bad.cell.def}.
\item If a subinterval $I^{(j)}_k$ is detected as an $\epsilon$-badcell, then bisect this subinterval and add the end points to the partition $P^{b}_{N_{j+1}}$, referred as the {\it badcells partition}.
\item Define the new partition $P_{N_{j+1}} = P_{N_{j}} \cup P^{b}_{N_{j+1}}$ and denote the subintervals of these partition as $I^{(j+1)}_k$, for $k=0,1,2,\cdots, N_{j+1}-1.$
\item Let $l_{*}: = \min\{|I^{(j+1)}_k|~|~k=0,1,2,\cdots, N_{j+1}-1\}$. If $l_*<\tau,$ then stop the iteration. Otherwise, repeat the above process for those intervals $I^{(j)}_k$ which are identified as $\epsilon$-badcells.
\end{enumerate}
The final outcome of the above process is the required {\it adaptive partition}, which we denote by $P^{*}_{N_J}$.
\subsection{Degree Adaptation}\label{degree.adaptation.subs}
The adaptive partition proposed above is expected to give an efficient way of obtaining the PiPCT approximation of a function $f\in L^{1}_{\omega}[a,b]$ with isolated singularities. In the construction of the adaptive partition, we have fixed $n_p$ and $n_q$ to be equal and equal in all the subintervals of the partition.  The numerical results depicted in Figure \ref{fig:PiPCheandChe_discondeg}(c) suggests that a more accurate approximation may also be obtained by choosing higher degrees for the numerator and denominator polynomials in the  $\epsilon$-badcells of the partition  $P^{*}_{N_J}$.  However, the error bounds \eqref{eq:totalL1errPade1}-\eqref{eq:totalL1errPade2} clearly shows that this is not the case in the region where the function is sufficiently smooth. More precisely, we see from \eqref{eq:totalL1errPade1} that the upper bound decreases as $d$ increases with $d<n$, whereas the inequality  \eqref{eq:totalL1errPade2} shows the opposite behaviour of the upper bound when $d>n$. This suggests us to choose $m< n/2$ sufficiently small in the smooth regions.  

\begin{figure}
\centering
\includegraphics[height=5cm,width=16cm]{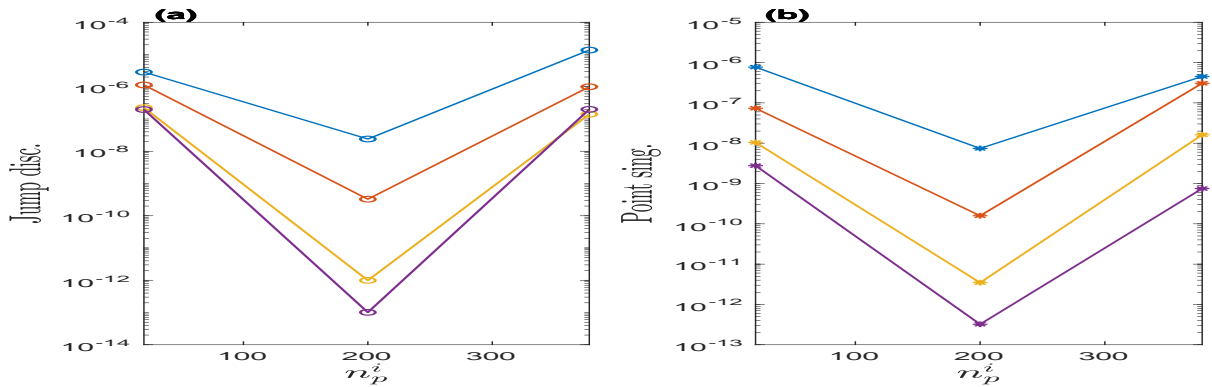}	
\caption{Pointwise error in PiPCT approximation of the function $f$ given by  \eqref{discondeg.eq} in a small neighborhood of \textbf{(a)} jump singularity at $x = -0.4$ and \textbf{(b)} point singularity at $x = 0.4.$ Different lines correspond to $N = 104, 208, 312, 416$, where $n_q^i=m=20$ and $n=200$. Three points lying on a graph (shown in symbol) corresponds to the case when (from left to right) $n_p = n_q,$ $n_p=n,$ and  $n_p=2n-n_q-1$, respectively.}
\label{fig:Adaptivetool}
\end{figure}

In this subsection, we propose a better choice of the degree of the numerator polynomial in the $\epsilon$-badcells. 
From Proposition \ref{prop:1}, we see that the coefficient vector $\bfm{q}$ of the denominator polynomial of $\mathsf{R}_{m,m}$ is the flip vector of the denominator polynomial of $\mathsf{R}_{2n-m-1,m}$  for a given $m<n-1.$  Since the coefficient vector $\bfm{p}$ of the numerator polynomial of a Pad\'e-Chebyshev type approximant is computed using $\bfm{q}$, we expect that the errors involved in approximating $f$ in a $\epsilon$-badcell by $\mathsf{R}_{m,m}$ and $\mathsf{R}_{2n-m-1,m}$ are almost the same. 
In particular, this result has been verified numerically for the function $f$ given by \eqref{discondeg.eq} when $m$ is sufficiently smaller than $n-1$ and is shown in Figure \ref{fig:Adaptivetool}. This figure suggests that the error attains its minimum when $n_p=n$ in this example.

Based on the above numerical observation, we suggest the following choice of the polynomial degrees in the adaptive partition obtained using the algorithm explained in subsection \ref{generation.adaptive.partition.ssec}.

\noindent{\it Choice of polynomial degrees in the adaptive algorithm:}\\
For a given number of quadrature points $n$, choose $m\ll n/2$.  Set $n_p=n$ and $n_q=m$ in all $\epsilon$-badcells of the adaptive partition $P^*_{N_J}$ and set $n_p=n_q=m$ for other cells of the partition. We call the resulting algorithm as {\it adaptive piecewise Pad\'e-Chebyshev type} (APiPCT) algorithm.
\begin{figure}[t]
\centering
\includegraphics[height=12.5cm,width=12cm]{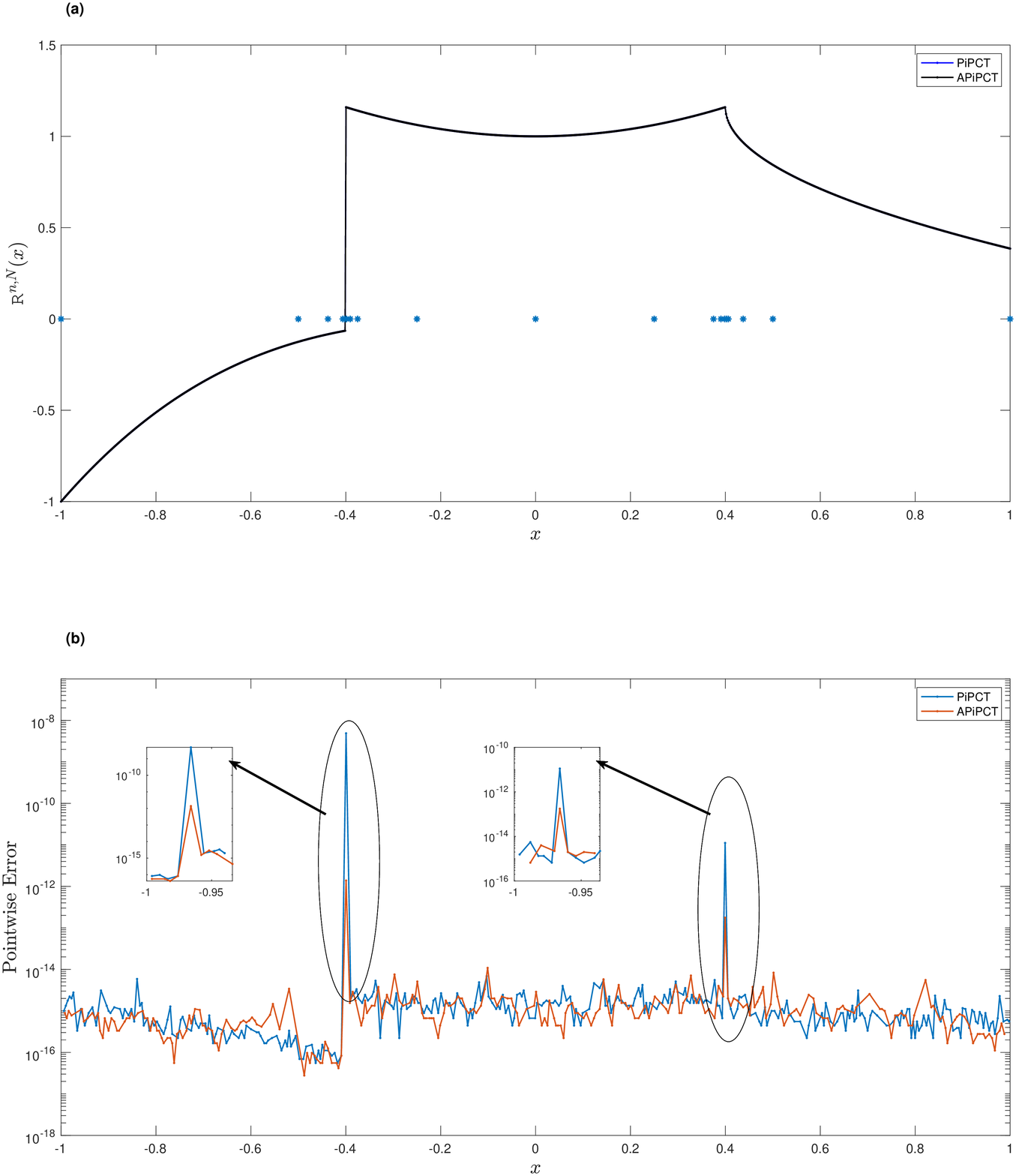}
\caption{\textbf{(a)} Depicts the comparison of the approximants obtained by the PiPCT algorithm and the APiPCT algorithm. \textbf{(b)} Depicts the comparison of the pointwise error of the approximants using the two algorithms for $N = 512,n = 100$, and for numerator and denominator degrees $n_p = n_q = 20$ for the function $f(t)$ given by \eqref{discondeg.eq}.}
\label{fig:plotErrPiPCAPiPC}
\end{figure}
\section{Numerical Experiment}\label{numerical.experiment.APiPCT.sec}
To validate the APiPCT algorithm, we compute the  PiPCT approximant of the function $f$ given by \eqref{discondeg.eq} using both  PiPCT and  APiPCT algorithms. In Figure \ref{fig:plotErrPiPCAPiPC} (a), we compare the approximants from both the algorithms.  Here, we choose $n=100$, $m=20,$ and $N=512$.  For the APiPCT algorithm, we choose $\epsilon =  10^{-2}$ and $\tau=(b-a)/N=1/256.$  With these parameters, the adaptive algorithm takes 18 cells as depicted by `*' symbol in Figure \ref{fig:plotErrPiPCAPiPC} (a).  Figure \ref{fig:plotErrPiPCAPiPC} (b) depicts the corresponding pointwise errors.  Here, we clearly observe that the error in the badcells are significantly reduced in the APiPCT algorithm when compared to that of PiPCT. This improvement is mainly because of the choice of numerator polynomial degree $n_p=n$ in the $\epsilon$-badcells.

We are also interested in studying the performance of the robust Pad\'e-Chebyshev  (RPCT) method of Gonnet {\it et al.} \cite{gon-etal_13a} in the adaptive piecewise algorithm developed in the above section. 
Observe that without any further modification, we can replace the PCT method by the RPCT method  in the adaptive algorithm. We denote the adaptive piecewise RPCT algorithm  by APiRPCT method. We observe that APiRPCT method also captures the singularities as sharply as APiPCT algorithm.

\begin{figure}[t]
	\centering
	\includegraphics[height=7cm,width=16cm]{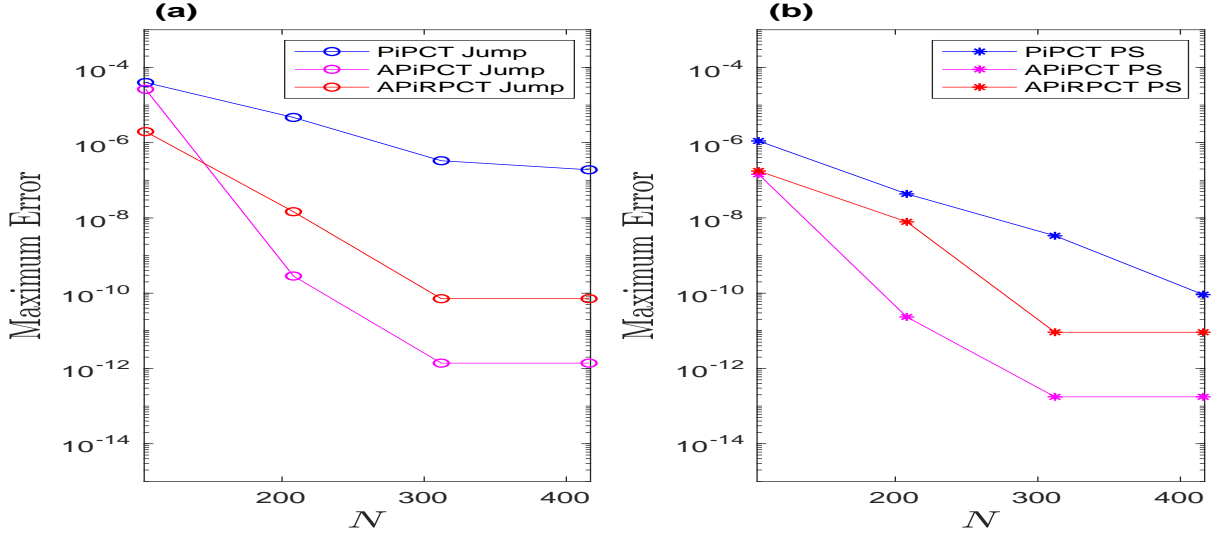}	
	\caption{Comparison between PiPCT, APiPCT and APiRPCT algorithms. \textbf{(a)} depicts the comparison in maximum error in a neighborhood of the jump discontinuity at $x = -0.4$ and \textbf{(b)} depicts the comparison in maximum error in a neighborhood of point singularity at $x = 0.4$,  for $N = 104, 208, 312, 416$, and for numerator and denominator degrees $n_p = n_q = 20$ for function $f(t)$ \eqref{discondeg.eq} on the domain $[-1,1]$. }
	\label{fig:Figure62}
\end{figure}

Figure \ref{fig:Figure62} \textbf{(a)} and \textbf{(b)} depicts the peaks of pointwise error in a vicinity of singularities present in the function \eqref{discondeg.eq} at $x = -0.4$ (jump discontinuity) and $x = 0.4$ (point singularity), respectively, for $N =
104, 208, 312, 416$. We can clearly see the performance of the three methods (PiPCT, APiPCT, and APiRPCT) near a singularity and the power of choosing degrees of numerator and denominator adaptively. Here, we understand that for a given $n_q$, setting $n_p=n$ in the adaptive methods decreases the maximum error in an $\epsilon$-badcell more rapidly than in PiPCT as $N$ increases. 

We also observe from Figure \ref{fig:Figure62} \textbf{(a)} and \textbf{(b)} that the maximum error in APiPCT decreases more rapidly than the APiRPCT method. This is because the RPCT method calculates the rank of the Toeplitz matrix (using singular value decomposition) and reduces the numerator and denominator degrees diagonally to reach to the final position (correspond to the minimum degree denominator) in upper left corner of a square block of the Pad\'e table (for block structure details see Gragg \cite{gra_72a} and Trefethen \cite{tre_84a}). Thus, in the process of minimising the occurrence of spurious pole-zero pair, the RPCT method decreases the denominator degree and therefore reduces the accuracy (as already been mentioned in Gonnet {\it et al.} \cite{gon-etal_13a}).  

\begin{figure}[t]
\centering
\includegraphics[height=7cm,width=10cm]{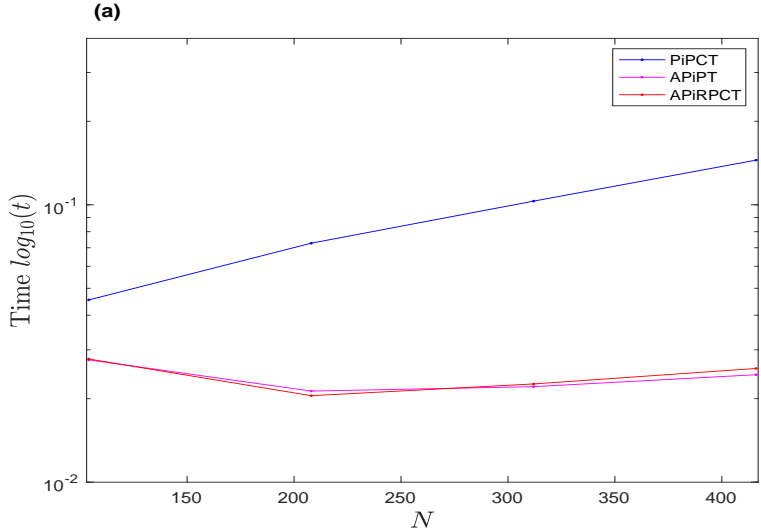}	
\caption{Comparison between the time taken by the PiPCT, the APiPCT and the APiRPCT algorithms to compute the approximants of the function $f(t)$ given by \eqref{discondeg.eq}, where $N = 104, 208, 312, 416$, and the numerator and the denominator degrees as $n_p = n_q = 20$. }
\label{fig:Efficiency}
\end{figure}

Finally, to study the efficiency of the adaptive piecewise algorithm, we compare the time taken by PiPCT, APiPCT, and APiRPCT algorithms to approximate the function \eqref{discondeg.eq}. Figure \ref{fig:Efficiency} depicts the comparison between the time taken by these three algorithms as we increase the number of partitions $N$ (correspondingly decreasing the tolerance parameter $\tau=(b-a)/N$). In this figure, we observe that the time taken by PiPCT approximation increases with $N$ while the time taken by the adaptive algorithms APiPCT and APiRPCT remain almost the same. Although the RPCT approximation does a repeated singular value decomposition to get a full rank matrix, finally the construction of the robust Pad\'e approximation is done possibly with a lower degree polynomials and hence the time taken in the repeated singular value decomposition is compensated in the Pad\'e construction.  Whereas, in APiPCT,  the Pad\'e approximation is done with higher order polynomials. Hence the time taken by the APiPCT and APiRPCT finally remained almost the same in this example.

\section{Comments on Froissart Doublets}\label{comment.FD.sec}
From theoretical point of view, it is known that (see for instance Baker and Peter \cite{bak-pet_82a}) a Pad\'e approximation accelerate the convergence of a truncated series. It can be used as a noise filter in signal processing. It significantly reduce the effect of Gibbs oscillations but not able to eliminate it completely. Apart from these (theoretical) properties of a Pad\'e (rational) approximation, there are difficulties in elucidating the approximation power of these approximants correctly. It is mainly because of the random occurrence of the  poles of the PCT approximant in the complex plane. In the case of a real valued function if these poles are sufficiently away from the unit circle then it may not effect the approximation. The main problem occurs when an approximant has a pole at the place where function has no singularities and in such cases one can not expect an accurate results (near the pole). These poles are referred as \textit{spurious poles}. Baker \textit{et al.} in \cite{bak-etal_61a} shows that spurious poles are isolated and always accompanied with a zero. The mentioned spurious pole-zero pair is also known as \textit{Froissart doublets}. The occurrence of Froissart doublets is a fundamental mathematical issue, these are the difficulties in establishing convergence theorems of Pad\'e approximants of order $[n_p/n_p]$, which cannot be achieved without a restriction of convergence in measure or capacity and not uniform convergence \cite{pom_73a,nut_70a}. On a computer in floating point arithmetic they even arise more often and known as \textit{numerical Froissart doublets} (see Nakatsukasa {\it et al.} \cite{nak-etal_18a}). 

There are at least two methods to recognize numerical Froissart doublets. In \cite{gon-etal_13a,gon-etal_11a,nak-etal_18a}, authors proposed a method to recognize those doublets by seeing the absolute value of their residual. Another method to recognize the spurious pole-zero pair is by seeing the distance between the two \cite{bes-per_09a}. In this section, we use the residual method to recognise the spurious poles.  

\begin{figure}[]
	\centering
	\includegraphics[height=5cm,width=5cm]{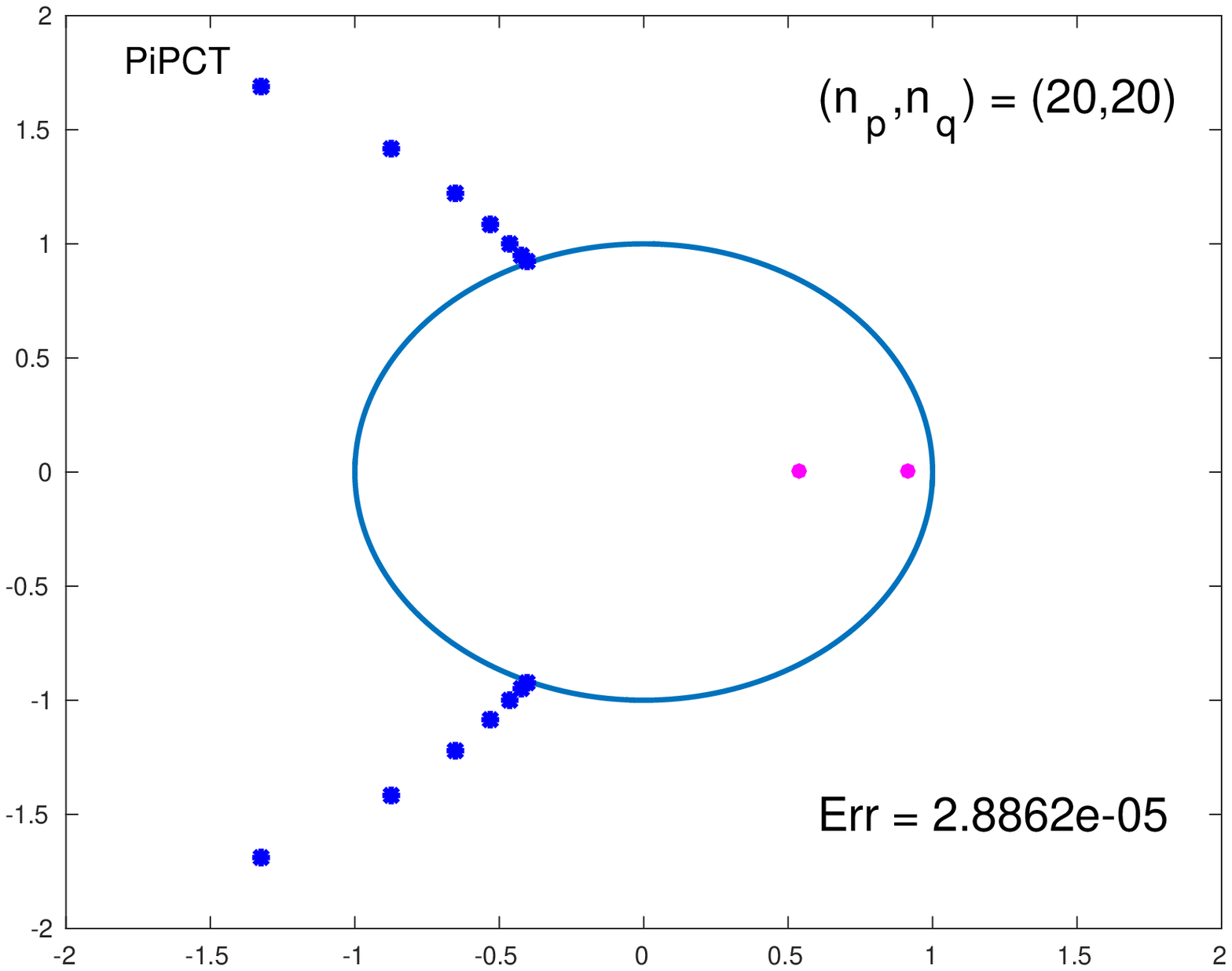}	
	\includegraphics[height=5cm,width=5cm]{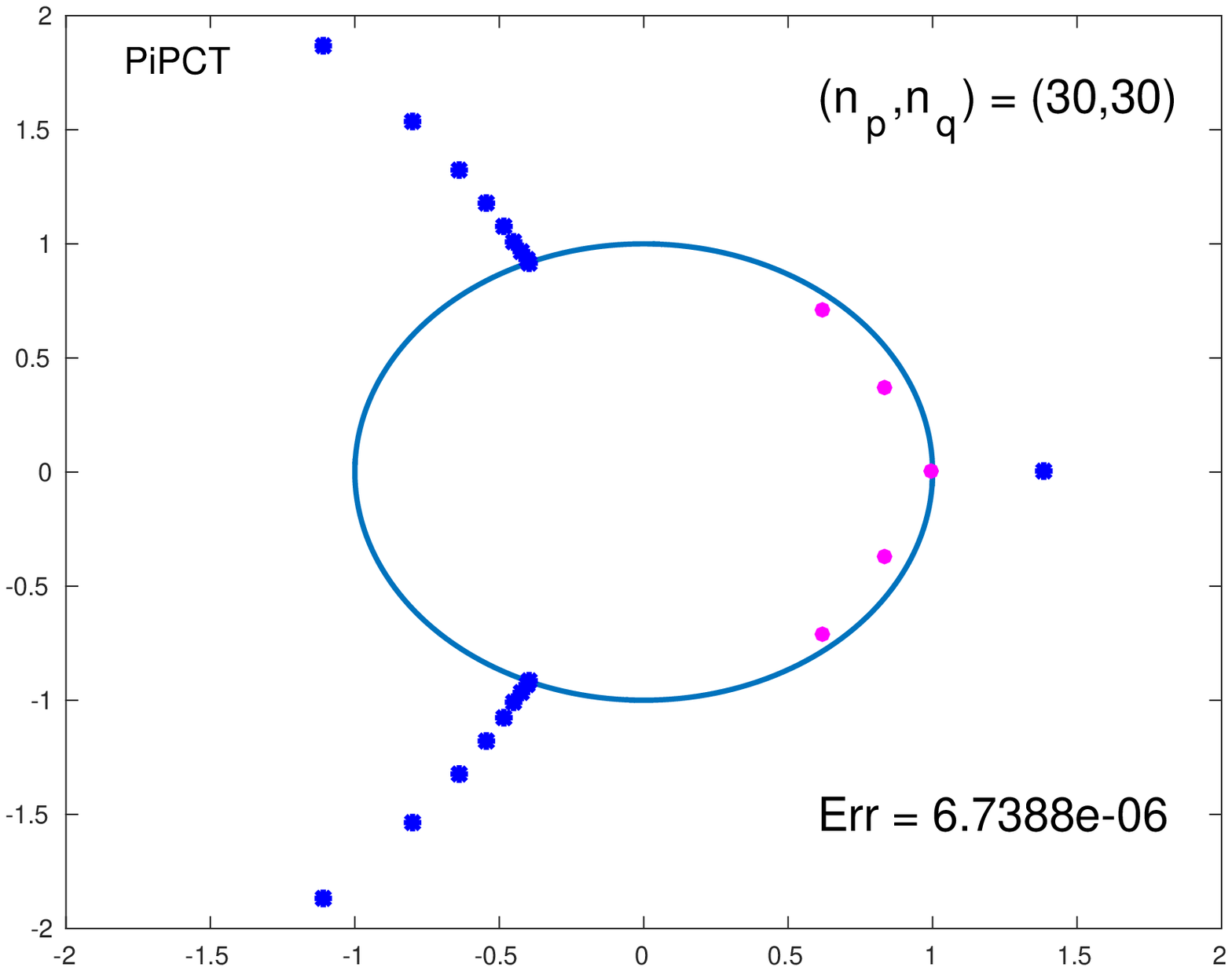}	
	\includegraphics[height=5cm,width=5cm]{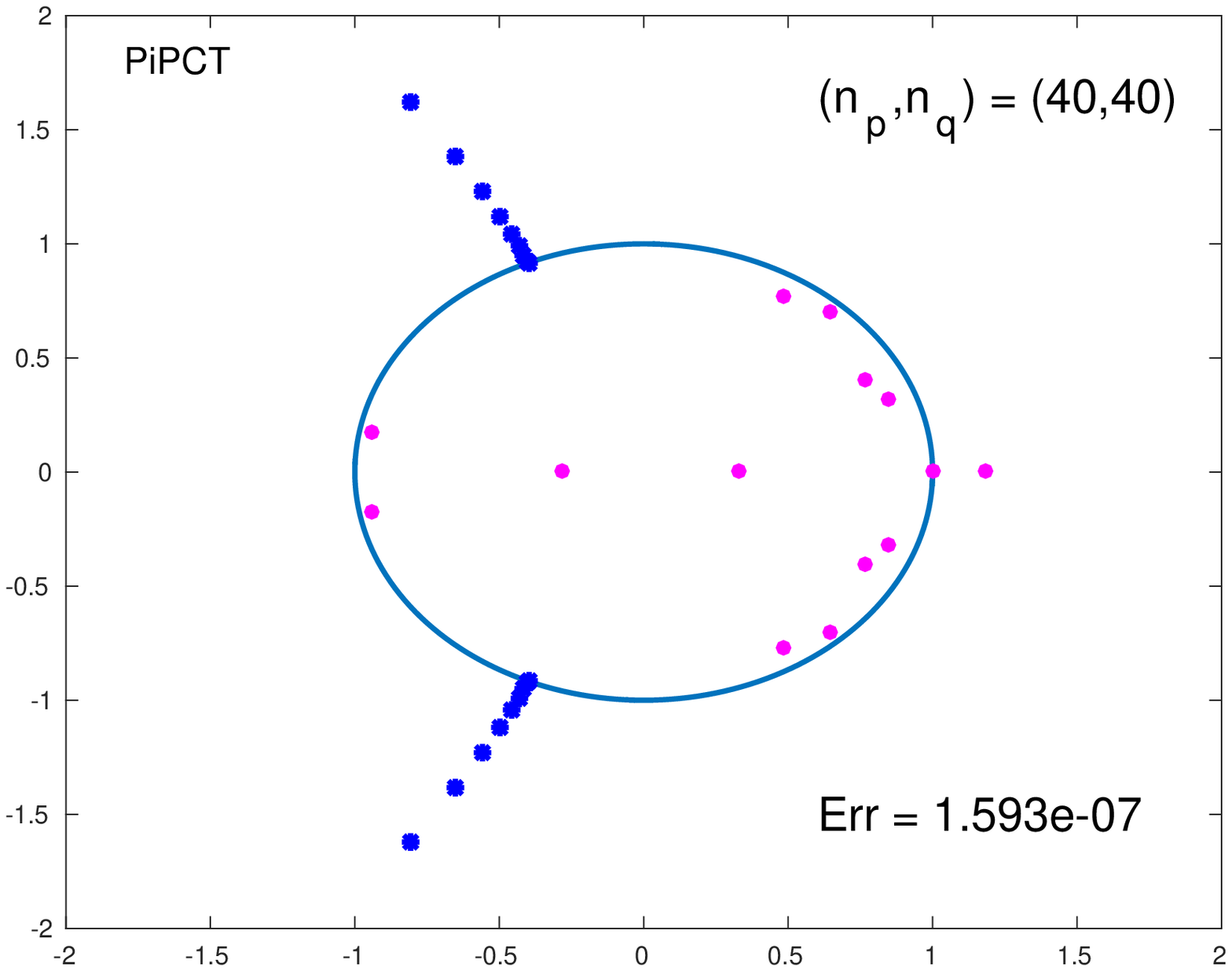}
	\includegraphics[height=5cm,width=5cm]{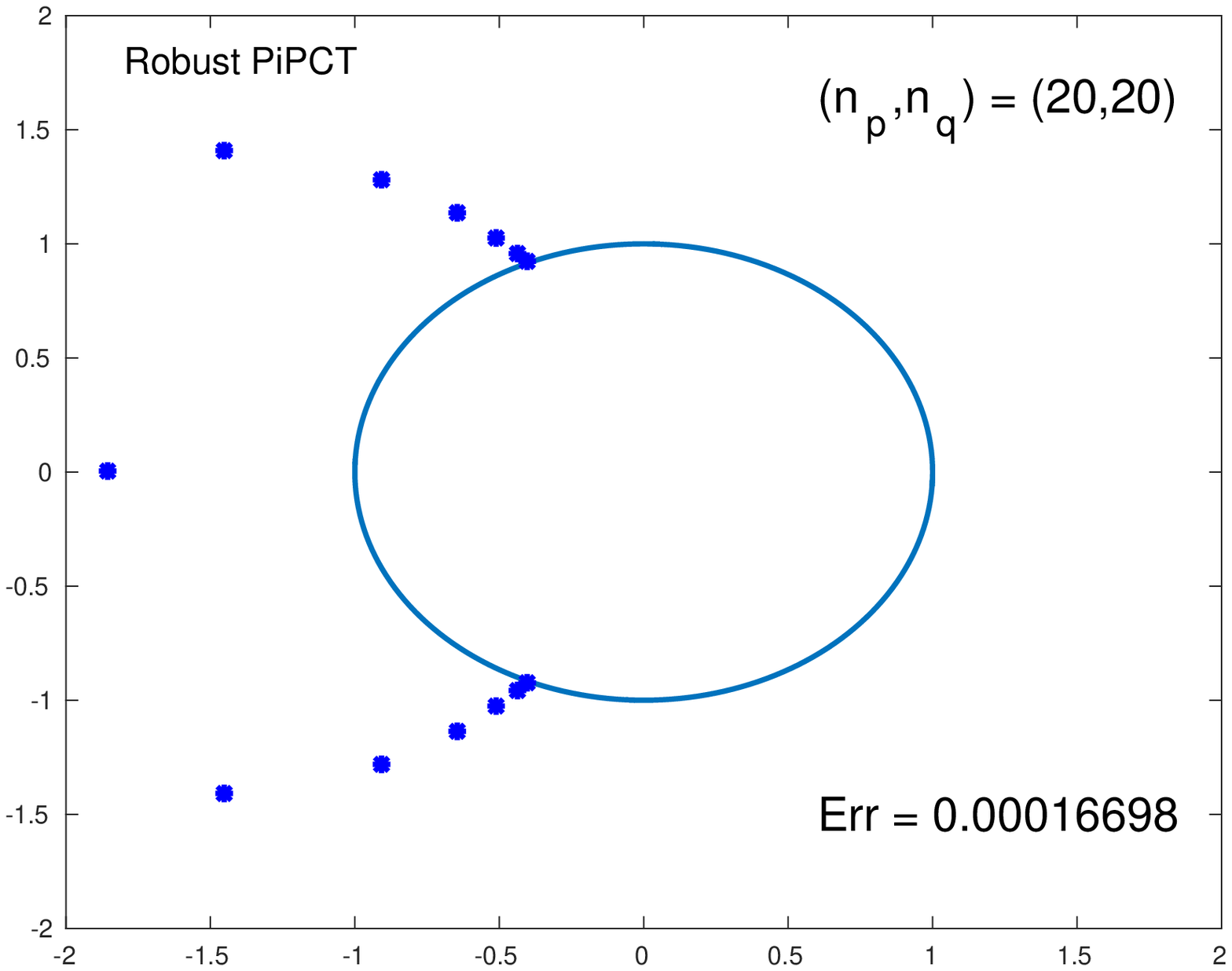}	
	\includegraphics[height=5cm,width=5cm]{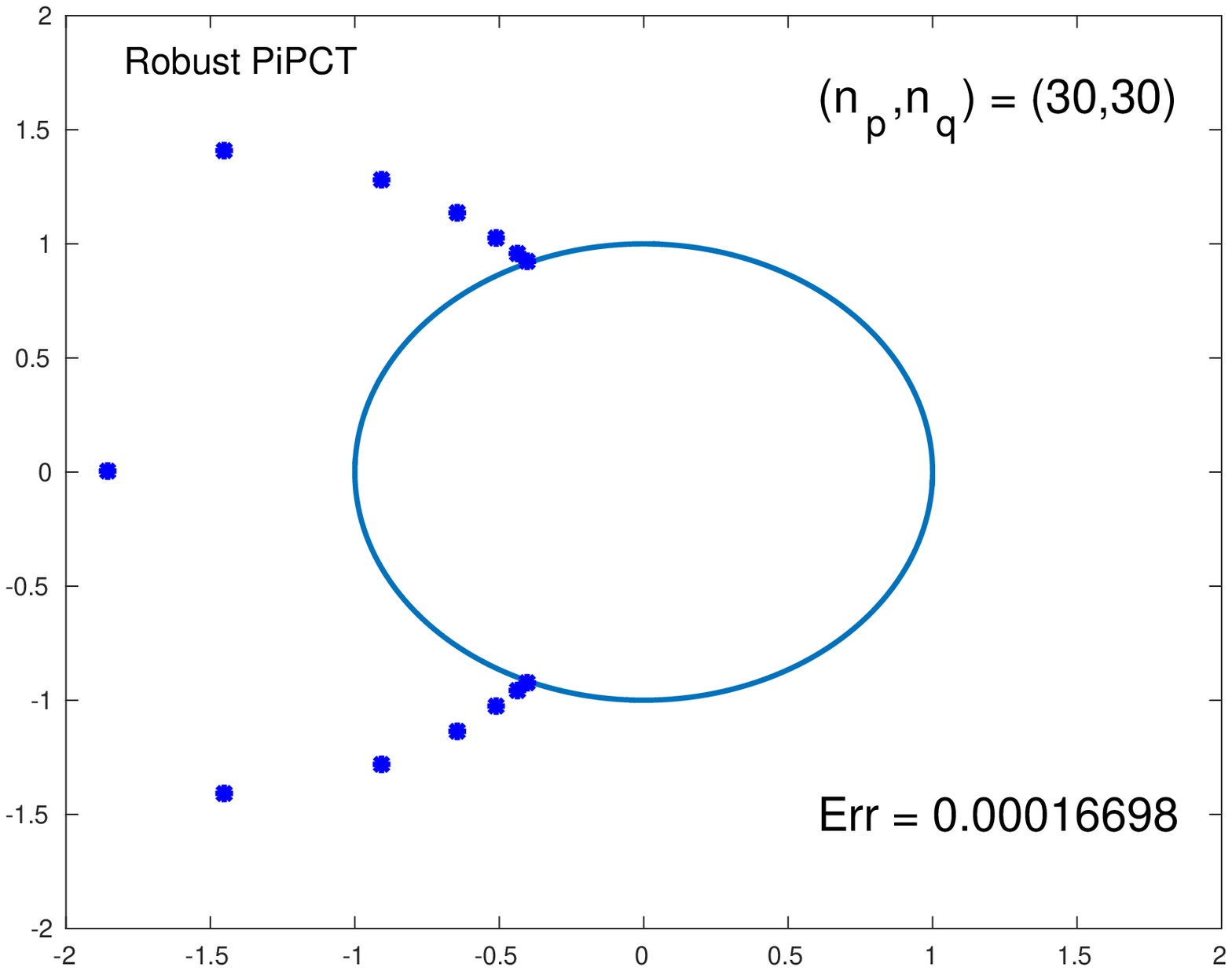}	
	\includegraphics[height=5cm,width=5cm]{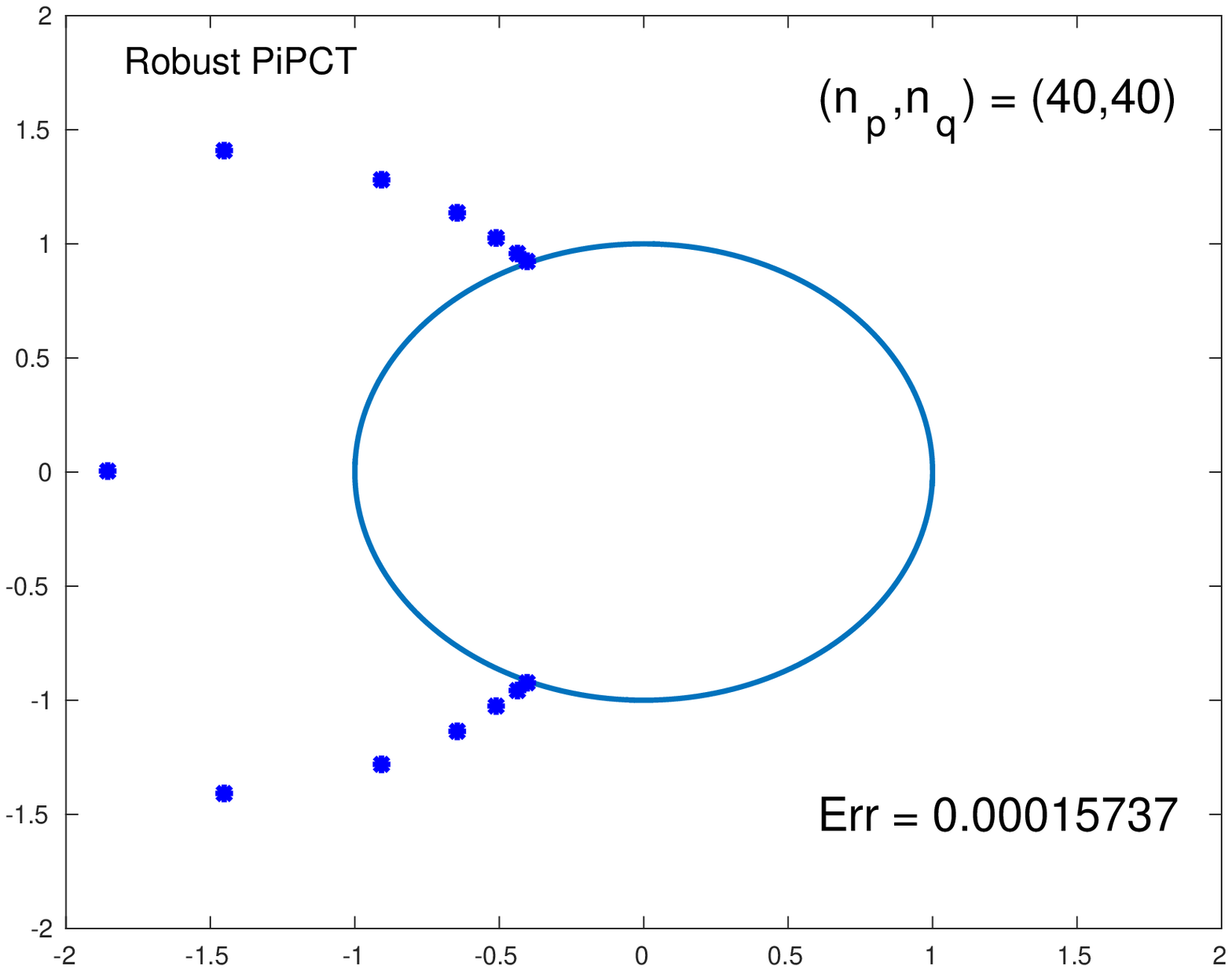}	
\caption{PiPCT and PiRPCT algorithms for $n_p = n_q = 20,30,40$ is performed. }
\label{fig:RPCTpoles}
\end{figure}

Figure \ref{fig:RPCTpoles} depicts the poles of PCT and RPCT approximations in $\epsilon$-badcell with jump discontinuity at $x = -0.4$. The figures in first row depicts the poles of the PCT approximation and the second row corresponds to the poles of the RPCT approximation.  Here pink circles denote the spurious poles and blue circles denote genuine poles.  From the figures of the first row, we observe that the spurious poles do occur in the  $\epsilon$-badcells, which shows that the  APiPCT is not free from Froissart doublets. We further observe that the number of spurious pole-zero pair increases as $n_p=n_q$ increases, along with the increase in the accuracy.   Also, we observe from second row of Figure \ref{fig:RPCTpoles} that the RPCT is free from spurious poles, with a compromise in accuracy.  This is clearly because the RPCT method eliminates spurious poles (in this example) by reducing the degrees of the numerator and the denominator polynomials.

\section{Conclusion} 
A piecewise Pad\'e-Chebyshev type (PiPCT) algorithm for approximating non-smooth functions on a bounded interval is proposed. A $L^1$-error estimate is obtained in the regions where the target function is smooth.  The error estimate shows that the order of converges depends on the degree of smoothness of the function as the number of partitions of the interval, $N\rightarrow \infty$. Numerical experiments are performed with a function involving both a jump discontinuity and a point singularity.  The convergence of the PiPCT method in the vicinity of singularities are demonstrated numerically both in the case of increasing $N$ values and increasing degrees of the numerator and the denominator polynomials of the Pad\'e-Chebyshev approximant. The numerical results clearly show the acceleration of the convergence of the approximant which is not the case with the global Pad\'e-Chebyshev type approximation. Another advantage of the PiPCT approach is that one does not need to know the location and the type of the singularities {\it a priori}, unlike in many methods that are proposed in the literature. The PiPCT algorithm is designed to work on a nonuniform mesh, which makes the algorithm more flexible for choosing a suitable adaptive partition and degrees.  A strategy based on the zeros of the denominator polynomial is used to identify a cell, called the {\it badcell}, where a possible singularity of the function lies.  Using this strategy, we have further developed an adaptive way of partitioning the interval so as to minimize the computational time without compromising the accuracy.  We also proposed a way to choose the numerator degree so as to improve the accuracy of the approximant in the bad cells and the resulting algorithm is called the adaptive piecewise Pad\'e-Chebyshev type (APiPCT) algorithm.  Numerical experiments are performed to validate the APiPCT algorithm in terms of accuracy and time efficiency.  A comparison study is performed where the  numerical results of PiPCT is compared with some recently developed similar methods like the singular Pad\'e-Chebyshev method and the robust Pad\'e-Chebyshev method.


\begin{thebibliography}{10}

\bibitem{ara-etal_05a}
F.~Arandiga, A.~Cohen, R.~Donat, and N.~Dyn.
\newblock Interpolation and approximation of piecewise smooth functions.
\newblock {\em SIAM J. Numer. Anal.}, 43(1):41--57, 2005.

\bibitem{bak-pet_82a}
G.~A. Baker and G.-M. Peter.
\newblock The convergence of sequences of pad{\'e} approximants.
\newblock {\em Journal of Mathematical Analysis and Applications}, 87(2):382 --
  394, 1982.

\bibitem{bak-etal_61a}
G.~A. Baker, J., J.~L. Gammel, and J.~G. Wills.
\newblock An investigation of the applicability of the {P}ad\'{e} approximant
  method.
\newblock {\em J. Math. Anal. Appl.}, 2:405--418, 1961.

\bibitem{bak_65a}
G.~A. Baker, Jr.
\newblock The theory and application of the {P}ad\'e approximant method.
\newblock In {\em Advances in {T}heoretical {P}hysics, {V}ol. 1}, pages 1--58.
  Academic Press, New York, 1965.

\bibitem{bak-etal_96a}
G.~A. Baker, Jr. and G.-M. Peter.
\newblock {\em Pad\'e approximants}, volume~59 of {\em Encyclopedia of
  Mathematics and its Applications}.
\newblock Cambridge University Press, Cambridge, second edition, 1996.

\bibitem{ban-gee_98a}
N.~S. Banerjee and J.~F. Geer.
\newblock Exponentially accurate approximations to periodic lipschitz functions
  based on fourier series partial sums.
\newblock {\em Journal of Scientific Computing}, 13(4):419--460, Dec 1998.

\bibitem{bar-etal_07a}
A.~Barkhudaryan, R.~Barkhudaryan, and A.~Poghosyan.
\newblock Asymptotic behavior of {E}ckhoff's method for {F}ourier series
  convergence acceleration.
\newblock {\em Anal. Theory Appl.}, 23(3):228--242, 2007.

\bibitem{bec-mat_15a}
B.~Beckermann and A.~C. Matos.
\newblock Algebraic properties of robust {P}ad\'{e} approximants.
\newblock {\em J. Approx. Theory}, 190:91--115, 2015.

\bibitem{bes-per_09a}
D.~Bessis and L.~Perotti.
\newblock Universal analytic properties of noise: introducing the {$J$}-matrix
  formalism.
\newblock {\em J. Phys. A}, 42(36):365202, 15, 2009.

\bibitem{cuy-wuy_87a}
A.~Cuyt and L.~Wuytack.
\newblock {\em Nonlinear methods in numerical analysis}, volume 136 of {\em
  North-Holland Mathematics Studies}.
\newblock North-Holland Publishing Co., Amsterdam, 1987.
\newblock Studies in Computational Mathematics, 1.

\bibitem{dev_98a}
R.~A. DeVore.
\newblock Nonlinear approximation.
\newblock In {\em Acta numerica, 1998}, volume~7 of {\em Acta Numer.}, pages
  51--150. Cambridge Univ. Press, Cambridge, 1998.

\bibitem{dev-kun_09a}
R.~A. DeVore and A.~Kunoth.
\newblock {\em Multiscale, Nonlinear and Adaptive Approximation: Dedicated to
  Wolfgang Dahmen on the Occasion of his 60th Birthday}.
\newblock Springer-Verlag Berlin Heidelberg, 1 edition, 2009.

\bibitem{dri-for_01a}
T.~A. Driscoll and B.~Fornberg.
\newblock A {P}ad\'e-based algorithm for overcoming the {G}ibbs phenomenon.
\newblock {\em Numer. Algorithms}, 26(1):77--92, 2001.

\bibitem{eck_95a}
K.~S. Eckhoff.
\newblock Accurate reconstructions of functions of finite regularity from
  truncated {F}ourier series expansions.
\newblock {\em Math. Comp.}, 64(210):671--690, 1995.

\bibitem{gee_95a}
J.~F. Geer.
\newblock Rational trigonometric approximations using fourier series partial
  sums.
\newblock {\em Journal of Scientific Computing}, 10(3):325--356, Sep 1995.

\bibitem{gon-etal_13a}
P.~Gonnet, S.~G\"{u}ttel, and L.~N. Trefethen.
\newblock Robust {P}ad\'{e} approximation via {SVD}.
\newblock {\em SIAM Rev.}, 55(1):101--117, 2013.

\bibitem{gon-etal_11a}
P.~Gonnet, R.~Pach\'{o}n, and L.~N. Trefethen.
\newblock Robust rational interpolation and least-squares.
\newblock {\em Electron. Trans. Numer. Anal.}, 38:146--167, 2011.

\bibitem{got-shu_95a}
D.~Gottlieb and C.-W. Shu.
\newblock On the {G}ibbs phenomenon. {IV}. {R}ecovering exponential accuracy in
  a subinterval from a {G}egenbauer partial sum of a piecewise analytic
  function.
\newblock {\em Math. Comp.}, 64(211):1081--1095, 1995.

\bibitem{got-shu_97a}
D.~Gottlieb and C-W. Shu.
\newblock On the {G}ibbs phenomenon and its resolution.
\newblock {\em SIAM Rev.}, 39(4):644--668, 1997.

\bibitem{gra_72a}
W.~B. Gragg.
\newblock The {P}ad\'{e} table and its relation to certain algorithms of
  numerical analysis.
\newblock {\em SIAM Rev.}, 14:1--16, 1972.

\bibitem{hes-etal_06a}
J.~S. Hesthaven, S.~M. Kaber, and L.~Lurati.
\newblock Pad\'e-{L}egendre interpolants for {G}ibbs reconstruction.
\newblock {\em J. Sci. Comput.}, 28(2-3):337--359, 2006.

\bibitem{kab-mad_05a}
S.~M. Kaber and Y.~Maday.
\newblock Analysis of some {P}ad\'e-{C}hebyshev approximants.
\newblock {\em SIAM J. Numer. Anal.}, 43(1):437--454, 2005.

\bibitem{kve_04a}
G.~Kvernadze.
\newblock Approximating the jump discontinuities of a function by its
  {F}ourier-{J}acobi coefficients.
\newblock {\em Math. Comp.}, 73(246):731--751, 2004.

\bibitem{lip-lev_10a}
Y.~Lipman and D.~Levin.
\newblock Approximating piecewise-smooth functions.
\newblock {\em IMA J. Numer. Anal.}, 30(4):1159--1183, 2010.

\bibitem{lit_03a}
G.~L. Litvinov.
\newblock Error autocorrection in rational approximation and interval
  estimates. [{A} survey of results].
\newblock {\em Cent. Eur. J. Math.}, 1(1):36--60, 2003.

\bibitem{maj_17a}
H.~Majidian.
\newblock On the decay rate of chebyshev coefficients.
\newblock {\em Applied Numerical Mathematics}, 113:44 -- 53, 2017.

\bibitem{mas-han_03a}
J.~C. Mason and D.~C. Handscomb.
\newblock {\em Chebyshev polynomials}.
\newblock Chapman \& Hall/CRC, Boca Raton, FL, 2003.

\bibitem{min-etal_07a}
M.~S. Min, S.~M. Kaber, and W.~S. Don.
\newblock Fourier-{P}ad\'e approximations and filtering for spectral
  simulations of an incompressible {B}oussinesq convection problem.
\newblock {\em Math. Comp.}, 76(259):1275--1290, 2007.

\bibitem{nak-etal_18a}
Y.~Nakatsukasa, O.~S\'ete, and L.~N. Trefethen.
\newblock The aaa algorithm for rational approximation.
\newblock {\em SIAM Journal on Scientific Computing}, 40(3):A1494--A1522, Jan
  2018.

\bibitem{nut_70a}
J.~Nuttall.
\newblock The convergence of pad{\'e} approximants of meromorphic functions.
\newblock {\em Journal of Mathematical Analysis and Applications}, 31(1):147 --
  153, 1970.

\bibitem{pom_73a}
C.~Pommerenke.
\newblock Pad\'{e} approximants and convergence in capacity.
\newblock {\em J. Math. Anal. Appl.}, 41:775--780, 1973.

\bibitem{riv_74a}
T.~J. Rivlin.
\newblock {\em The {C}hebyshev polynomials}.
\newblock Wiley-Interscience [John Wiley \& Sons], New York-London-Sydney,
  1974.
\newblock Pure and Applied Mathematics.

\bibitem{sta_98a}
H.~Stahl.
\newblock Spurious poles in pad{\'e} approximation.
\newblock {\em Journal of Computational and Applied Mathematics}, 99(1):511 --
  527, 1998.
\newblock Proceeding of the VIIIth Symposium on Orthogonal Polynomials and
  Thier Application.

\bibitem{tad_07a}
E.~Tadmor.
\newblock Filters, mollifiers and the computation of the {G}ibbs phenomenon.
\newblock {\em Acta Numer.}, 16:305--378, 2007.

\bibitem{tam-etal_12a}
A.~L. Tampos, J.~E.~C. Lope, and J.~S. Hesthaven.
\newblock Accurate reconstruction of discontinuous functions using the singular
  {P}ad\'e-{C}hebyshev method.
\newblock {\em IAENG Int. J. Appl. Math.}, 42(4):242--249, 2012.

\bibitem{tre_84a}
L.~N. Trefethen.
\newblock Square blocks and equioscillation in the pad{\'e}, walsh, and cf
  tables.
\newblock In Peter~Russell Graves-Morris, Edward~B. Saff, and Richard~S. Varga,
  editors, {\em Rational Approximation and Interpolation}, pages 170--181,
  Berlin, Heidelberg, 1984. Springer Berlin Heidelberg.

\bibitem{xia-etal_10a}
S.~Xiang, X.~Chen, and H.~Wang.
\newblock Error bounds for approximation in {C}hebyshev points.
\newblock {\em Numer. Math.}, 116(3):463--491, 2010.

\end{thebibliography}
\end{document}